\documentclass[preprint,12pt]{elsarticle}

\usepackage{graphics}
\usepackage{hyperref}
\usepackage{amsmath}
\usepackage{amssymb}
\usepackage{color}
\usepackage{listings}
\usepackage{textcomp}
\definecolor{listinggray}{gray}{0.9}
\definecolor{lbcolor}{rgb}{0.9,0.9,0.9}
\newcommand{\go}{\stackrel{\circ }{\mathfrak{g}}}

\newcommand{\co}[1]{\stackrel{\circ }{#1}}

\newcommand{\gf}{\mathfrak{g}}
\newcommand{\af}{\mathfrak{a}}
\newcommand{\bff}{\mathfrak{b}}
\newcommand{\afb}{\mathfrak{a}_{\bot}}
\newcommand{\hf}{\mathfrak{h}}
\newcommand{\hfg}{\hf_{\gf}}

\newcommand{\pf}{\mathfrak{p}}
\newcommand{\aft}{\widetilde{\mathfrak{a}}}

\newcounter{bla}

\journal{Computer Physics Communications}

\begin{document}

\begin{frontmatter}



\title{{\bf Affine.m} -- {\it Mathematica} package for computations in representation theory of finite-dimensional and affine Lie algebras}


\author[a,b]{Anton Nazarov\corref{author}}
%\author[a,b]{Second Author}
%\author[b]{Third Author}

\cortext[author] {Corresponding author.\\\textit{E-mail address:} antonnaz@gmail.com}
\address[a]{Department of High Energy Physics, Faculty of physics, SPb State University\\ 198904, Sankt-Petersburg, Russia}
\address[b]{Chebyshev Laboratory, Faculty of Mathematics and Mechanics, SPb State University\\ 199178, Saint-Petersburg, Russia}

\begin{abstract}
In this paper we present {\bf Affine.m} --- a program for computations in representation theory of finite-dimensional and affine Lie algebras and describe implemented algorithms. The algorithms are based on the properties of weights and Weyl symmetry.  Computation of weight multiplicities in irreducible and Verma modules, branching of representations and tensor product decomposition are the most important problems for us. These problems have numerous applications in physics and we provide some examples of these applications. The program is implemented in the popular computer algebra system {\it Mathematica} and works with finite-dimensional and affine Lie algebras. 

% A submitted program is expected to be of benefit to other physicists or physical chemists, or be an exemplar of good programming practice, or illustrate new or novel programming techniques which are of importance to some branch of computational physics or physical chemistry.
%
% Acceptable program descriptions can take different forms. The following Long Write-Up structure is a suggested structure but it is not obligatory. Actual structure will depend on the length of the program, the extent to which the algorithms or software have already been described in literature, and the detail provided in the user manual.

%Your manuscript and figure sources should be submitted through the Elsevier Editorial System (EES) by using the online submission tool at \\
% http://www.ees.elsevier.com/cpc.

%In addition to the manuscript you must supply: the program source code; job control scripts, where applicable; a README file giving the names and a brief description of all the files that make up the package and clear instructions on the installation and execution of the program; sample input and output data for at least one comprehensive test run; and, where appropriate, a user manual. These should be sent, via email as a compressed archive file, to the CPC Program Librarian at cpc@qub.ac.uk.

\end{abstract}


\end{frontmatter}

% \linenumbers

% Computer program descriptions should contain the following
% PROGRAM SUMMARY.

{\bf PROGRAM SUMMARY}%/NEW VERSION PROGRAM SUMMARY}
  %Delete as appropriate.

\begin{small}
\noindent
{\em Manuscript Title:}{\bf Affine.m} -- {\it Mathematica} package for computations in representation theory of finite-dimensional and affine Lie algebras                                       \\
{\em Authors:}Anton Nazarov                                                \\
{\em Program Title:}Affine.m                                          \\
{\em Catalogue identifier:}  AENA\textunderscore v1\textunderscore 0                                 \\
  %Leave blank, supplied by Elsevier.
{\em Licensing provisions:}Standard CPC licence, \url{http://cpc.cs.qub.ac.uk/licence/licence.html}                                   \\
{\em No. of lines in distributed program, including test data, etc.:} 24 844\\ 
{\em No. of bytes in distributed program, including test data, etc.:} 1 045 908\\
{\em Distribution format:} tar.gz\\
{\em Programming language:}Mathematica                                   \\
{\em Computer:}i386-i686, x86\textunderscore 64                                               \\
  %Computer(s) for which program has been designed.
{\em Operating system:} Linux, Windows, Mac OS, Solaris                                       \\
  %Operating system(s) for which program has been designed.
{\em RAM:} 5-500 Mb                                              \\
  %RAM in bytes required to execute program with typical data.
%{\em Number of processors used:}                              \\
  %If more than one processor.
%{\em Supplementary material:}                                 \\
  % Fill in if necessary, otherwise leave out.
  % Please give some freely chosen keywords that we can use in a
  % cumulative keyword index.
{\em Classification:} 5 Computer Algebra, 4.2 Other algebras and groups                                         \\
  %Classify using CPC Program Library Subject Index, see (
  % http://cpc.cs.qub.ac.uk/subjectIndex/SUBJECT_index.html)
  %e.g. 4.4 Feynman diagrams, 5 Computer Algebra.
%{\em External routines/libraries:}                            \\
  % Fill in if necessary, otherwise leave out.
%{\em Subprograms used:}                                       \\
  %Fill in if necessary, otherwise leave out.
%{\em Catalogue identifier of previous version:}*              \\
  %Only required for a New Version summary, otherwise leave out.
%{\em Journal reference of previous version:}*                  \\
  %Only required for a New Version summary, otherwise leave out.
%{\em Does the new version supersede the previous version?:}*   \\
  %Only required for a New Version summary, otherwise leave out.
{\em Nature of problem:}\\
  %Describe the nature of the problem here.
Representation theory of finite-dimensional Lie algebras has many applications in different branches of physics, including elementary particle physics, molecular physics, nuclear physics. Representations of affine Lie algebras appear in string theories and two-dimensional conformal field theory used for the description of critical phenomena in two-dimensional systems. Also Lie symmetries play major role in a study of quantum integrable systems.
   \\
{\em Solution method:}\\
  %Describe the method solution here.
We work with weights and roots of finite-dimensional and affine Lie algebras and use Weyl symmetry extensively. Central problems which are the computations of weight multiplicities, branching and fusion coefficients are solved using one general recurrent algorithm based on generalization of Weyl character formula. We also offer alternative implementation based on the Freudenthal multiplicity formula which can be faster in some cases.
   \\
%{\em Reasons for the new version:}*\\
  %Only required for a New Version summary, otherwise leave out.
%   \\
%{\em Summary of revisions:}*\\
  %Only required for a New Version summary, otherwise leave out.
%   \\
{\em Restrictions:}\\
  %Describe any restrictions on the complexity of the problem here.
Computational complexity grows fast  with the rank of an algebra, so computations for algebras of ranks greater than 8 are not practical.
   \\
{\em Unusual features:}\\
  %Describe any unusual features of the program/problem here.
We offer the possibility to use a traditional mathematical notation for the objects in representation theory of Lie algebras in computations if {\bf Affine.m} is used in the {\it Mathematica} notebook interface.
   \\
%{\em Additional comments:}\\
  %Provide any additional comments here.
%   \\
{\em Running time:}\\
  %Give an indication of the typical running time here.
From seconds to days depending on the rank of the algebra and the complexity of the representation.
   \\

% \begin{thebibliography}{0}
% \bibitem{1}Reference 1         % This list should only contain those items referenced in the
% \bibitem{2}Reference 2         % Program Summary section.
% \bibitem{3}Reference 3         % Type references in text as [1], [2], etc.
%                                % This list is different from the bibliography at the end of
%                                % the Long Write-Up.
% \end{thebibliography}
% * Items marked with an asterisk are only required for new versions
% of programs previously published in the CPC Program Library.\\
\end{small}

\section{Introduction}
\label{intro}

Representation theory of Lie algebras is of central importance for different areas of physics and mathematics. In physics Lie algebras are used to describe symmetries of quantum and classical systems. Computational methods in representation theory have a long history \cite{belinfante1989survey}, there exist numerous software packages for computations related to Lie algebras \cite{simplie}, \cite{vanleeuwen1994lsp}, \cite{stembridge1995mps,coxweyl}, \cite{fischbacher2002ilp}, \cite{Fuchs:1996dd}.

Most popular programs \cite{simplie}, \cite{vanleeuwen1994lsp}, \cite{fischbacher2002ilp}, \cite{coxweyl} are created to study representation theory of simple finite-dimensional Lie algebras. The main computational problems are the following:
\begin{enumerate}
\item Construction of a root system which determines the properties of Lie algebra including its commutation relations.
\item Weyl group traversal which is important due to Weyl symmetry of root system and characters of representations.
\item Calculation of weight multiplicities, branching and fusion coefficients, which are essential for construction and study of representations.
\end{enumerate}
There are well-known algorithms for these tasks \cite{moody1982fast}, \cite{stembridge2001computational}, \cite{belinfante1989survey}, \cite{casselman1994machine}.
The third problem is the most computation intensive. There are two different recurrent algorithms which are based on the Weyl character formula and the Freudenthal multiplicity formula. In this paper we analyze them both.

Infinite-dimensional Lie algebras also have a growing number of applications in physics for example in conformal field theory and in a study of  integrable systems. But infinite-dimensional algebras are much harder to investigate and the number of available computer programs is much smaller in this case.

Affine Lie algebras \cite{kac1990idl} constitute an important and tractable class of infinite-dimensional Lie algebras. They are constructed as central extensions of loop algebras of (semi-simple) finite-dimensional Lie algebras and appear naturally in a study of Wess-Zumino-Witten and coset models in conformal field theory \cite{Walton:1999xc}, \cite{difrancesco1997cft}, \cite{Goddard198588}, \cite{Dunbar:1992gh}.

The structure of affine Lie algebras allows us to adapt for them the computational algorithms created for finite-dimensional Lie algebras  \cite{Fuchs:1996dd}, \cite{gannon2001algorithms}, \cite{kass1990ala}. The book \cite{kass1990ala} with the tables of multiplicities and other computed characteristics of affine Lie algebras and representations was published in 1990. But we are not aware of software packages for popular computer algebra systems which can be used to extend these results.
We address this issue and present {\bf Affine.m} -- a {\it Mathematica} package for computations in representation theory of affine and finite-dimensional Lie algebras.  We describe the features and limitations of the package in the present paper.  We also provide a representation-theoretical background of the implemented algorithms and present the examples of computations relevant to physical applications.

The paper starts with an overview of Lie algebras and their representation theory (Sec. \ref{sec:theor-backgr}). Then we describe the datastructures of {\bf Affine.m} used to present different objects related to Lie algebras and representations (Sec. \ref{sec:core-datastructures}) and discuss the implemented algorithms (Sec. \ref{sec:comp-algor}). The next section consists of physically interesting examples (Sec. \ref{sec:examples}). The paper is concluded with the discussion of possible extensions and refinements (Sec. \ref{sec:conclusion}).

\section{Theoretical background}
\label{sec:theor-backgr}

In this section we recall necessary definitions and present formulae used in computations. 

\subsection{Lie algebras of finite and affine types}
\label{sec:lie-algebras-finite}

A {\it Lie algebra} $\gf$ is a vector space with a bilinear operation $[\cdot,\cdot]:\gf\otimes\gf\to \gf$, which is called a {\it Lie bracket}. If we choose a basis $X_{i}$ in $\gf$ we can specify the commutation relations by the {\it structure constants} $C_{ijk}$:
\begin{equation}
  \label{eq:1}
  [X^{i},X^{j}]=\sum_{k} C^{ij}_{k} X^{k}
\end{equation}
A Lie algebra is {\it simple} if it contains no non-trivial ideals with respect to a commutator. A {\it semisimple} Lie algebra is a direct sum of simple Lie algebras. In the present paper we treat simple and semisimple Lie algebras. 

A {\it Cartan subalgebra}  $\hfg$ is a nilpotent subalgebra of  $\gf$ that coincides with its normalizer. 
We denote the elements of a basis of $\hfg$ by $H^{i}$.

The Killing form on $\gf$ gives a non-degenerate bilinear form $(\cdot,\cdot)$ on $\hfg$ which can be used to identify $\hfg$ with the subspace of the dual space $\hfg^{*}$ of linear functionals on $\hfg$. {\it Weights} are the elements of $\hfg^{*}$ and are denoted by Greek letters $\mu,\nu, \omega, \lambda\dots$

A special choice of a basis gives a compact description of the commutation relations (\ref{eq:1}). This basis can be encoded by the root system which is discussed in Section \ref{sec:weights-roots} (See also \cite{humphreys1997introduction,humphreys1992reflection}).

The {\it loop algebra} $L\gf=\gf\otimes \mathbb{C}[t,t^{-1}]$, corresponding to semisimple Lie algebra $\gf$, has commutation relations
\begin{equation}
  \label{eq:6}
  [X^{i}t^{n},X^{j}t^{m}]=t^{n_+m}\sum_{k}C^{ij}_{k}X^{k}
\end{equation}
The central extension leads to the appearance of an additional term
\begin{equation}
  \label{eq:7}
   [X^{i}t^{n}+\alpha c,X^{j}t^{m}+\beta c]=t^{n+m}\sum_{k}C^{ij}_{k}X^{k}+(X^{i},X^{j})n\delta_{n+m,0}c
\end{equation}
This algebra $\hat\gf=\gf\otimes\mathbb{C}[t,t^{-1}]\oplus\mathbb{C}c$ is called a non-twisted {\it affine Lie algebra} \cite{kac1990idl}, \cite{wakimoto2001idl,wakimoto2001lectures}, \cite{kass1990ala}.

We do not treat twisted affine Lie algebras in the present paper.

\subsection{Modules, weights and roots}
\label{sec:weights-roots}

Let $\gf$ be a finite-dimensional or an affine Lie algebra.

Then the $\gf$-module is a vector space $V$ together with a bilinear map $\gf \times V\to V$ such that
\begin{equation}
  \label{eq:2}
  [x,y]\cdot v = x\cdot(y\cdot v) - y\cdot(x\cdot v), \quad \mbox{for}\; x,y\in \gf, v\in V
\end{equation}
A representation of an algebra $\gf$ on a vector space $V$ is a homomorphism $\gf\to gl(V)$ from $\gf$ to a Lie algebra of endomorphisms of the  vector space $V$ with the commutator as the bracket.

For an arbitrary representation it is possible to diagonalize the operators corresponding to Cartan generators $H^{i}$ simultaneously by a special choice of basis $\{v_{j}\}$ in $V$:
\begin{equation}
  \label{eq:3}
  H^{i}\cdot v_{j}=\nu_{j}^{i}v_{j}
\end{equation}
The eigenvalues $\nu^{i}_{j}$ of Cartan generators on an element of basis $v_{j}$ determine a weight $\nu_{j}\in \hfg^{*}$ such that $\nu_{j}(H^{i})=\nu_{j}^{i}$. A vector $v\in V$ is called the weight vector of the weight $\lambda$ if $H v=\lambda_{j}(H)v,\; \forall H\in \hf$. The  weight subspace consists of all weight vectors $V_{\lambda}=\{v\in V: H v=\lambda_{j}(H)v,\; \forall H\in \hf\}$. The weight multiplicity $m_{\lambda}=\mathrm{mult}(\lambda)=\mathrm{dim} V_{\lambda}$ is the dimension of the weight subspace.

The structure of a module is determined by the set of weights since the action of generators $E^{\alpha}$ on weight vectors is
\begin{equation}
  \label{eq:5}
  E^{\alpha}\cdot v_{\lambda} \propto v_{\lambda+\alpha}
\end{equation}
The module structure can be encoded by the formal character of the module
\begin{equation}
  \label{eq:10}
  \mathrm{ch}V=\sum_{\lambda}m_{\lambda} e^{\lambda}
\end{equation}
The character  $\mathrm{ch}V\in \mathcal{E}$ is an element of algebra $\mathcal{E}$ generated by formal exponents of weights.
The character can be specialized by taking its value on some element $\xi$ of $\hf$.

Any Lie algebra is its own module with respect to a special kind of representation. The action that
defines this representation is called {\it adjoint} and is given by the bracket $ad_{X} Y=[X,Y]$. 
{\it Roots} are weights of the adjoint representation of $\gf$. They encode the commutation relations of algebra in the following way. Denote by $\Delta$ the set of roots. For each $\alpha\in \Delta$ there exist a root $-\alpha\in \Delta$ and the generators $E^{\alpha}, E^{-\alpha}$ such that
\begin{align}
  \label{eq:4}
  &  [H^{i},E^{\alpha}]=\alpha^{i}E^{\alpha} \\
  &\left[E^{\alpha},E^{\beta}\right]=
  \begin{cases}
    N_{\alpha,\beta} E^{\alpha+\beta}, & \mbox{if}\; \alpha+\beta\in \Delta\\
    \frac{2}{(\alpha,\alpha)} \sum_{i}\alpha^{i} H^{i},&  \mbox{if}\; \alpha=-\beta\\
    0,&\mbox{otherwise}
  \end{cases}
\end{align}

Given the root system $\Delta$ we can choose the set of positive roots. This is a subset  $\Delta^{+}\subset \Delta$ such that for each root $\alpha\in\Delta$ exactly one of the roots $\alpha, -\alpha$ is contained in $\Delta^{+}$ and for any two distinct positive roots $\alpha, \beta\in \Delta^{+}$ such that $\alpha+\beta\in \Delta$ their sum is also positive $\alpha+\beta\in\Delta^{+}$.
The elements of $-\Delta^{+}$ are called negative roots.

A positive root is {\it simple} if it cannot be written as a sum of positive roots. The set of simple roots $\Phi=\left\{\alpha_{i}\right\}$ is a basis in $\hfg^{*}$ and each root can be written as $\alpha=\sum_{i}n_{i}\alpha_{i}$ with all $n_{i}$ non-negative or non-positive simultaneously. In the case of a finite-dimensional Lie algebra $\gf$ simple roots are numbered from 1 to the rank of the algebra $i=1,\dots,r,\quad r=\mathrm{rank}(\gf)$. Enumerating simple roots with an index $i$ we introduce lexicographic ordering in the root system $\Delta$. The highest root with respect to this ordering is denoted by  $\theta=\sum_{i=1,\dots,r} a_i \alpha_i$, the coefficients $a_i$ are called {\it marks}. $\theta$ is also the highest weight  of the adjoint module (See section \ref{sec:high-weight-modul}). {\it Comarks} are the numbers equal to $a_i^{\vee}=\frac{(\alpha_i,\alpha_i)}{2} a_i$.

Although for an affine Lie algebra $\hat\gf$ the set of roots $\Delta$ is infinite the set of simple roots $\Phi$ is finite and its elements are denoted by $\alpha_{0},\dots \alpha_{r}$ where $r=\mathrm{rank}(\gf)$. The roots $\alpha_1,\dots, \alpha_r$ are the roots of the underlying finite-dimensional Lie algebra $\go$. The root $\alpha_0=\delta-\theta$ is the difference of the {\it imaginary root} $\delta$ and the highest root  $\theta$ of the algebra $\go$.
Note that root multiplicity $\mathrm{mult}(\alpha)$ for an affine Lie algebra root can be greater than one.

A subalgebra  $\bff_{+}\subset \gf$ spanned by the generators $H^{i}, E^{\alpha}$ for positive roots $\alpha\in \Delta^{+}$ is called a {\it Borel subalgebra}.

A {\it parabolic subalgebra}  $\pf_{I}\supset \bff_{+}$ contains a Borel subalgebra and is  generated by some subset of simple roots $\{\alpha_{j}:j\in I, I\subset \{1\dots r\}\}$. It is spanned by the subset of generators $\{H^{i}\}\cup \{E^{\alpha}:\alpha\in \Delta^{+}\}\cup \{E^{-\alpha}: \alpha\in\Delta^{+}, \alpha=\sum_{j\in I} n_{j} \alpha_{j}\}$.
A {\it regular subalgebra} $\af\subset\gf$ is determined by the root system $\Delta_{\af}$ with the set of simple roots $\{\beta_{i}, i=1,\dots,r_{\af}\}$ being a subset of the set of roots $\{\alpha_{1},\dots,\alpha_{r}\}\cup \{-\theta\}$ .

The {\it Weyl group} $W_{\gf}$ is generated by reflections $\{s_{i}:\hfg^{*}\to\hfg^{*}\}$ corresponding to simple roots $\{\alpha_{i}\}$:
\begin{equation}
  \label{eq:8}
  s_{i}\cdot\lambda=\lambda-\frac{2(\alpha_{i},\lambda)}{(\alpha_{i},\alpha_{i})}\alpha_{i}
\end{equation}
The root system and characters of representation are invariant with respect to the Weyl group  action. A root system can be reconstructed from the set of simple roots by the Weyl group transformations.

Weyl groups are finite for finite-dimensional Lie algebras and finitely generated for affine Lie algebras.

Consider an action of the element $s_{\alpha}s_{\alpha+\delta}$ of the Weyl group of affine Lie algebra  $\hat\gf$ for $\alpha$ being a simple root of the underlying finite-dimensional Lie algebra $\gf$. Using definition \eqref{eq:8} it is easy to see that $s_{\alpha}s_{\alpha+\delta} \cdot \lambda=\lambda+\frac{2}{(\alpha,\alpha)}\alpha+\left(\frac{(\alpha,\alpha)}{2 (\lambda,\delta)}+(\lambda,\alpha)\right) \delta$. So the Weyl group $W_{\hat\gf}$ can be presented as a semidirect product of the Weyl group $W_{\gf}$ of $\gf$ and the set of translations corresponding to the roots of $\gf$. 

A Weyl group element can be presented as a product of elementary reflections in multiple ways. The number of elementary reflections in the shortest sequence representing an element $w\in W_{\gf}$ is called the {\it length} of $w$ and is denoted by $l(w)$. We also use the notation $\epsilon(w)=(-1)^{l(w)}$ for the parity of the number of Weyl reflections generating $w$. 

A fundamental domain $\bar{C}$ for the Weyl group $W_{\gf}$ action  on $\hfg^{*}$ is determined by the requirement $\xi\in \bar{C}\Leftrightarrow (\xi,\alpha_{i})\geq 0$ for all simple roots $\alpha_{i}$. It is called the {\it main Weyl chamber}. 

A {\it Cartan matrix} $A$ is defined by products of simple roots
\begin{equation}
  \label{eq:9}
  A_{ij}=\frac{2(\alpha_{i},\alpha_{j})}{(\alpha_{j},\alpha_{j})}
\end{equation}
and can be used for a compact description of Lie algebra commutation relations in the Chevalley basis \cite{humphreys1997introduction}, \cite{fulton1991representation}, \cite{bourbaki2002lie}.

The form \eqref{eq:9} induces a basis dual  to the simple
roots basis. It is called the {\it fundamental weights basis}. We denote
its elements by $\omega_i$:
\begin{equation}
  \label{eq:20}
  \langle\omega_i,\alpha_j\rangle=\frac{2(\omega_{i},\alpha_{j})}{(\alpha_{j},\alpha_{j})}=\delta_{ij}
\end{equation}
For a finite-dimensional Lie algebra there are $r$ fundamental weights, $i=1,\dots, r$. For an affine Lie algebra we have an additional fundamental weight $\omega_0=\lambda$, $(\lambda,\delta)=1, \; (\lambda,\lambda)=(\delta,\delta)=0$. Other fundamental weights are equal to $\omega_i=a_i^v \omega_0 +\co{\omega_i}$, where $\co{\omega_i}$ is a fundamental weight of the finite-dimensional Lie algebra $\gf$.

 The sum of fundamental weights $\rho=\sum_{i} \omega_{i}$ is called a {\it Weyl vector}. It is an important tool in representation theory.

\subsection{Highest weight modules}
\label{sec:high-weight-modul}

We consider finitely generated $\gf$-modules $V$ such that $V=\bigoplus_{\xi\in \hfg^{*}} V_{\xi}$, where each $V_{\xi}$ is finite dimensional and there exists a finite set of weights $\lambda_{1},\dots \lambda_{s}$ which generates the weight system of $V$, i.e. if $\mathrm{dim}V_{\xi}\neq 0$ then $\xi=\lambda_{i}-\sum_{k=1,\dots, r} n_{k}\alpha_{k}$ where $n_{k}\in \mathbb{Z}_{+}$ (See \cite{humphreys2008representations}, \cite{carter2005lie}).

A highest weight module $V^{\mu}$ contains a single highest weight $\mu$,  the other weights are obtained by subtractions of linear combinations of simple roots $\lambda=\mu-n_{1}\alpha_{1}-\dots-n_{r}\alpha_{r},\; n_{k}\in \mathbb{Z} _{+}$.

The simplest type of highest weight modules is the Verma module
$M^{\mu}$. Its space can be defined as a module
\begin{equation}
  \label{eq:17}
  M^{\mu}=U(\gf)\underset{U(\bff_{+})}{\otimes} D^{\mu}(\bff_{+}),
\end{equation}
with respect to a multiplication in $U(\gf)$ and
$\underset{U(\bff_{+})}{\otimes}$ means that the action of elements of $U(\bff_{+})$ ``falls through'' the left part of the tensor product onto the right part. Here $\bff_{+}$ is the Borel subalgebra, $D^{\mu}(\bff_{+})$ is a representation of $\bff_{+}$ such that $D(E^{\alpha})=0,\; D(H)=\mu(H)$ for any positive root $\alpha$.
Elements of $\gf$ act from the left and we should commute all the elements of $\bff_{+}$ to the right, so that they can act on the space $D^{\lambda}(\bff_{+})$.

Weight multiplicities in Verma modules can be found using the Weyl
character formula
\begin{equation}
  \label{eq:11}
  \mathrm{ch} M^{\mu}=\frac{e^{\mu}}{\prod_{\alpha\in \Delta^{+}} \left( 1-e^{-\alpha}\right)^{\mathrm{mult}(\alpha)}}=\frac{e^{\mu}}{\sum_{w\in W} \epsilon(w) e^{w\rho-\rho}}
\end{equation}
Here we have used the Weyl denominator identity
\begin{equation}
  \label{eq:12}
  R:=\prod_{\alpha\in \Delta^{+}} \left( 1-e^{-\alpha}\right)^{\mathrm{mult}(\alpha)}=\sum_{w\in W} \epsilon(w) e^{w\rho-\rho},
\end{equation}
and $\epsilon \left( w\right) :=\det \left( w\right)$ is equal to 
the parity of the sequence of Weyl reflections generating $w$.

A Verma module $M^{\mu}$ has the unique maximal submodule and the
unique nontrivial simple quotient $L^{\mu}$ which is an
{\it irreducible highest weight module}. 

Irreducible highest weight modules have no non-trivial submodules. 
The Weyl character formula for an irreducible highest weight module $L^{\mu}$ is
\begin{equation}
  \label{eq:13}
  \mathrm{ch} L^{\mu}=\frac{\sum_{w\in W} \epsilon(w) e^{w(\mu+\rho)-\rho}}{\sum_{w\in W}\epsilon(w) e^{w\rho-\rho}}=\sum_{w\in W} \epsilon(w)\; \mathrm{ch} M^{w(\mu+\rho)-\rho}
\end{equation}
Thus the character of an irreducible highest weight module can be
seen as a combination of characters of Verma modules. ( This
fact is a consequence of the Bernstein-Gelfand-Gelfand resolution
(\cite{bernstein1976category,bernstein1971structure}, see also
\cite{humphreys2008representations}).)

Construction of a generalized Verma module is analogous to (\ref{eq:17}), but the representation of the Borel subalgebra is substituted by a representation of a parabolic subalgebra $\pf_{I}\supset \bff_{+}$ generated by some subset $\{\alpha_{I}\}$ of simple roots $I\subset \{1,\dots, r\}$:
\begin{equation*}
M_{I}^{\mu}=U\left( \gf\right)\otimes _{U\left( \pf_{I}\right) }L_{\pf_{I}}^{\mu}.
\end{equation*}
Introduce a formal element $R_{I}:=\prod_{\alpha \in \Delta
^{+}\setminus \Delta _{\pf_{I}}^{+}}\left( 1-e^{-\alpha }\right)
^{\mathrm{mult}(\alpha )}$. Then the character of a generalized
Verma module can be written as
\begin{equation}
  \label{eq:18}
  \mathrm{ch}M_{I}^{\mu}=\frac{1}{R_{I}}\mathrm{ch}L_{\pf_{I}}^{\mu }.
\end{equation}

%% Define external border of irreducible representation here

We can use the Weyl character formula to obtain recurrent relations for weight multiplicities -- important tools for calculations \cite{il2010folded,kulish4sfa}. 

For irreducible highest-weight modules the recurrent relation has the following form
\begin{equation}
\label{eq:14}
m_{\xi }=-\sum_{w\in W\setminus e}\epsilon (w)m_{\xi
-\left( w(\rho )-\rho \right) }+\sum_{w\in W}\epsilon
(w)\delta _{\left( w(\mu +\rho )-\rho \right) ,\xi }.
\end{equation}
Formulae for Verma and generalized Verma modules differ only in the second term on the right-hand side. In the case of Verma module it is just $\delta_{\xi,\mu}$. For a generalized Verma module the summation in the second term on the right-hand side of \eqref{eq:14} is over the Weyl subgroup generated by the reflections corresponding to the roots $\{\alpha_{I}\}$.

Another recurrent formula can be obtained from a study of Casimir
element action on irreducible highest weight modules
\cite{humphreys1997introduction}:
\begin{equation}
  \label{eq:15}
  m_{\lambda}=\frac{2}{(\mu+\rho)^{2}-(\lambda+\rho)^{2}}\sum_{\alpha\in \Delta^{+}}\sum_{k\geq 1} (\lambda+k\alpha,\alpha)m_{\lambda+k\alpha}.
\end{equation}
It is called the Freudenthal multiplicity formula.
Note that it is applicable only to irreducible modules.

We discuss the use of formulae \eqref{eq:14} and \eqref{eq:15} for the computations in section \ref{sec:comp-algor}. 

Now consider an algebra $\gf$ and a reductive subalgebra
$\af\subset \gf$. Simple roots $\beta_{i}$ of the subalgebra $\af$
can be presented as linear combinations of $\gf$-algebra roots
$\alpha_{j}$: $\beta_{i}=\sum_{j=1,\dots,r_{\gf}}k_{j}
\alpha_{j},\ j=1,\dots,r_{\af}$.

Each irreducible $\gf$-module is also an $\af$-module, although
$L^{\mu}_{\gf}$ is in general not irreducible as an $\af$-module. 
It can be
decomposed into a direct sum of irreducible $\af$-modules:
\begin{equation}
  \label{eq:16}
  L^{\mu}_{\gf}=\bigoplus_{\nu}b^{\mu}_{\nu}L^{\nu}_{\af}
\end{equation}
The coefficients this decomposition are called branching coefficients. 

It is possible to calculate branching coefficients by constructing and successively  subtracting the submodules $L^{\nu}_{\af}$. 
This traditional approach has serious limitations especially in case of affine Lie algebras. We discuss them in the end of section \ref{sec:comp-algor}. 

Now we describe an alternative approach which is based on recurrent properties of branching coefficients. 
But before we proceed to these recurrent relations we need several additional definitions.

For a subalgebra $\af\subset \gf$ we introduce the subalgebra
$\afb$. Consider a root subspace $\hf_{\perp \af}^{\ast }$
orthogonal to $\hf_{\af}$,
\begin{equation*}
\hf_{\perp \af}^{\ast }:=\left\{ \eta \in \hf^{\ast }|\forall
h\in \hf_{\af};\eta \left( h\right) =0\right\} ,
\end{equation*}
and the roots (correspondingly -- positive roots) of $\gf$ orthogonal
to the roots of $\af$,
\begin{eqnarray}
\Delta _{\afb} &:&=\left\{ \beta \in \Delta _{\gf}|\forall
h\in \hf_{\af};\beta \left( h\right) =0\right\} ,
\label{delta a ort} \\
\Delta _{\afb}^{+} &:&=\left\{ \beta ^{+}\in \Delta _{\gf%
}^{+}|\forall h\in \hf_{\af};\beta ^{+}\left( h\right) =0\right\} .
\notag
\end{eqnarray}
Let $W_{\afb}$ be a subgroup of $W$ generated by
the reflections $w_{\beta }$ with the roots $\beta \in \Delta _{\af_{\perp
}}^{+}$. The subsystem $\Delta _{\afb}$ determines a 
subalgebra $\afb$ with the Cartan subalgebra $\hf_{\af%
_{\perp }}$.

The Cartan subalgebra $\frak{h}$ can be decomposed in the following way:  $\frak{h}=\frak{\frak{h}_{\af}}\oplus
\frak{h}_{\afb}\oplus \frak{h}_{\perp }$

We also introduce the notations
\begin{eqnarray}
\widetilde{\frak{a}_{\perp }} :=\frak{a}_{\perp }\oplus \frak{h}_{\perp }
\qquad
\widetilde{\frak{a}} :=\frak{a}\oplus \frak{h}_{\perp }.
\end{eqnarray}

For $\af$ and $\afb$ we consider the
corresponding Weyl vectors, $\rho _{\af}$ and $\rho _{\af_{\perp
}} $ and compose the so called ''defects'' $\mathcal{D}_{\af}$ and $\mathcal{%
D}_{\afb}$ of the injection:
\begin{equation}
\mathcal{D}_{\af}:=\rho _{\af}-\pi _{\af}\rho , \qquad
\mathcal{D}_{\afb}:=\rho _{\afb}-\pi _{\af%
_{\perp }}\rho .  \label{defect-ort}
\end{equation}

%%!!! Define P^+
For $\mu \in P^{+}$ consider the linked weights $\left\{ \left(
w(\mu +\rho )-\rho \right) |w\in W\right\} $ and their projections
to
$h_{\afb}^{\ast }$ additionally shifted by the defect $-%
\mathcal{D}_{\afb}$:
\begin{equation*}
\mu _{\afb}\left( w\right) :=\pi _{\afb}\left[
w(\mu +\rho )-\rho \right] -\mathcal{D}_{\afb},\quad w\in W.
\end{equation*}
Among the weights $\left\{ \mu _{\af_{\perp
}}\left( w\right) |w\in W\right\} $ one can always choose those located in
the fundamental chamber $\overline{C_{\afb}}$. Let $U$ be the
set of representatives $u$ for the classes $W/W_{\afb}$ such
that

\begin{equation}
U:=\left\{ u\in W|\quad \mu _{\afb}\left( u\right) \in
\overline{C_{\afb}}\right\} \quad .  \label{U-def}
\end{equation}
Thus we can form the subsets:
\begin{equation}
\mu _{\widetilde{\mathfrak{a}}}\left( u\right) :=\pi _{\widetilde{%
\mathfrak{a}}}\left[ u(\mu +\rho )-\rho \right] +\mathcal{D}_{\af%
_{\perp }},\quad u\in U,  \label{mu-a}
\end{equation}
and
\begin{equation}
\mu _{\afb}\left( u\right) :=\pi _{\afb}\left[
u(\mu +\rho )-\rho \right] -\mathcal{D}_{\afb},\quad u\in U.
\label{mu-a-tilda}
\end{equation}

Notice that the subalgebra $\mathfrak{a}_{\bot}$ is regular by definition
since it is built on a subset of roots of the algebra $\mathfrak{g}$.

Denote by  $k_{\xi }^{\left( \mu \right) }$ signed branching coefficients defined as follows. If $\xi\in \bar C_{\af}$ is in the main Weyl chamber $k_{\xi}^{(\mu)}=b^{(\mu)}_{\xi}$, otherwise $k_{\xi}^{(\mu)}=\epsilon(w) b^{(\mu)}_{w (\xi+\rho_{\af})-\rho_{\af}}$ where $w\in W_{\af}$ is such that $w (\xi+\rho_{\af})-\rho_{\af}\in \bar C_{\af}$. 

Now we can use the Weyl character formula to write a
recurrent relation \cite{2010arXiv1007.0318L} for the signed branching
coefficients $k_{\xi }^{\left( \mu \right) }$ corresponding to an
injection $\af\hookrightarrow \gf$:
\begin{equation}
\begin{array}{c}
k_{\xi }^{\left( \mu \right) }=-\frac{1}{s\left( \gamma _{0}\right) }\left(
\sum_{u\in U}\epsilon (u)\;\dim \left( L_{\afb}^{\mu _{\af%
_{\perp }}\left( u\right) }\right) \delta _{\xi -\gamma _{0},\pi _{%
\widetilde{\af}}(u(\mu +\rho )-\rho )}+\right.  \\
\left. +\sum_{\gamma \in \Gamma _{\widetilde{\af}\rightarrow \gf%
}}s\left( \gamma +\gamma _{0}\right) k_{\xi +\gamma }^{\left( \mu \right)
}\right) .
\end{array}
\label{recurrent-rel}
\end{equation}
The recursion is governed by the set $\Gamma _{\af\rightarrow \gf}$ called the injection fan. The latter is defined by the
carrier set $\left\{ \xi \right\} _{\af\rightarrow \gf}$ for the
coefficient function $s(\xi )$
\begin{equation*}
\left\{ \xi \right\} _{\widetilde{\af}\rightarrow \gf}:=\left\{
\xi \in P_{\widetilde{\af}}|s(\xi )\neq 0\right\}
\end{equation*}
appearing in the expansion
\begin{equation}
\prod_{\alpha \in \Delta ^{+}\setminus \Delta _{\bot }^{+}}\left( 1-e^{-\pi
_{\widetilde{\af}}\alpha }\right) ^{\mathrm{mult}(\alpha )-\mathrm{mult}%
_{\af}(\pi _{\widetilde{\af}}\alpha )}=-\sum_{\gamma \in P_{%
\widetilde{\af}}}s(\gamma )e^{-\gamma };\quad
\end{equation}
The weights in $\left\{ \xi \right\} _{\widetilde{\af}\rightarrow \gf}$ are to be shifted by $\gamma _{0}$ -- the lowest vector in $\left\{ \xi
\right\} $ -- and the zero element is to be eliminated:
\begin{equation}
\Gamma _{\af\rightarrow \gf}=\left\{ \xi -\gamma
_{0}|\xi \in \left\{ \xi \right\} \right\} \setminus \left\{ 0\right\} .
\end{equation}
The formula (\ref{eq:14}) is a particular case of the recurrent relation for branching coefficients (\ref{recurrent-rel}) in the case of a Cartan subalgebra $\af=\hfg$.

If the root system of $\afb$ is generated by some subset of $\gf$
simple roots $\alpha_{1},\dots,\alpha_{r}$ then the recurrent
relation (\ref{recurrent-rel}) is connected with the generalized
Bernstein-Gelfand-Gelfand resolution for parabolic Verma modules
\cite{2011arXiv1102.1702L}.

Another particular case of this formula is connected with tensor
product decompositions. Consider the tensor product of two
irreducible $\gf$-modules $L^{\mu}\otimes L^{\nu}$. It is also a
$\gf$-module but not irreducible in general. So
\begin{equation}
  \label{eq:19}
  L^{\mu}\otimes L^{\nu}=\bigoplus_{\gamma} f^{\mu\nu}_{\gamma}L^{\gamma}
\end{equation}
The coefficients $f^{\mu\nu}_{\gamma}$ are called the fusion
coefficients. The problem of computation of fusion coefficients is
equivalent to a branching problem for the diagonal subalgebra
$\gf\subset \gf\oplus \gf$ (see \cite{LyakhovskyPostnova2011}). So our implementation of a recurrent algorithm can be used to decompose tensor products (See Section \ref{sec:tens-prod-decomp}).

In the case of affine Lie algebras $\gf, \af$  the multiplicities
$m_{\nu}$ and the branching coefficients $b^{(\mu)}_{\nu}$
can be regarded as the coefficients in the power series decomposition
of string and branching functions correspondingly:
\begin{align}
  \label{eq:21}
  &\sigma_{\nu}(q)=\sum_{n=0}^{\infty} m_{\nu-n\delta} q^n, \quad \nu=\sum_j c_j \omega_j,\quad c_j\geq 0\\
  & b_{\nu}(q)=\sum_{n=0}^{\infty} b_{\nu-n\delta} q^n,\quad  \nu=\sum_j c_j \omega_j, \quad c_j\geq 0
\end{align}
String and branching functions have  modular and analytic properties which are important for conformal field theory, especially in  coset models and the study of CFT on higher genus surfaces \cite{kac1988modular}, \cite{difrancesco1997cft}, \cite{Walton:1999xc}, \cite{walton1989conformal}.

{\bf Affine.m}  calculates weight multiplicities and branching coefficients for affine Lie algebras up to some finite grade. We present examples of computations in Sections \ref{sec:string-funct-affine}, \ref{sec:branch-funct-coset}. Now we proceed to the description of the datastructures and the algorithms implemented in {\bf Affine.m}.

\section{Core datastructures}
\label{sec:core-datastructures} Having introduced necessary
mathematical objects, problems and relations we now describe the
related datastructures of {\bf Affine.m}. Although {\it
Mathematica} is untyped language it is possible to create
structured objects and do the type checks with patterns
\cite{shifrinmathematica}, \cite{maeder2000computer}.
\subsection{Weights}
\label{sec:weights}

Weights are represented  by two datastructures: \lstinline{finiteWeight} for finite-dimensional Lie algebras and \lstinline{affineWeight}  for affine.

Internally the finite weight is a \lstinline{List} with the
\lstinline{Head} \lstinline{finiteWeight}, its components are the
coordinates of the weight in the orthogonal Bourbaki  basis
\cite{bourbaki2002lie}.

An affine weight is an extension of a finite weight by supplying it
with the level and grade coordinates. There is a set of functions
defined for finite and affine weights. The complete list can be
found in the  online help of the package. The most important
are the definitions of an addition, a multiplication by a number and a
scalar product (bilinear form) for weights. These definitions
allow us to use a traditional notation with {\bf Affine.m}:
\begin{lstlisting}
  w=makeFiniteWeight[{1,0,3}];
  v=makeFiniteWeight[{3,2,1}];
  2*w+v==makeFiniteWeight[{5,2,7}]
  w.v==6
\end{lstlisting}

The use of an orthogonal basis in the internal structure of weights allows us to work with weights without a complete specification of a root system which is useful for a study of branching, since the roots of the subalgebra can be specified by hand.

\subsection{Root systems}
\label{sec:root-systems}

To specify an algebra of finite or affine type it is enough to fix
its root system. Root systems are represented by two datatypes
\lstinline{finiteRootSystem} and \lstinline{affineRootSystem}. The
latter is an extension of the former. We offer several different
constructors for these datastructures. It is possible to specify
the set of simple roots explicitly, for example to study the
subalgebra $B_2\subset B_4$ we can use the definition
\begin{lstlisting}
  b2b4=makeFiniteRootSystem[ { {1,-1,0,0}, {0,1,0,0} } ]
\end{lstlisting}
There are constructors for the root systems of simple finite-dimensional Lie algebras:
\begin{lstlisting}
  b2=makeSimpleRootSystem[B,2]
\end{lstlisting}
We use the typographic features of the {\it Mathematica} frontend to offer a mathematical notation for simple Lie algebras:

\begin{lstlisting}[mathescape=true]
  $B_2$ == makeFiniteRootSystem[ { {1, -1}, {0, 1} }]
\end{lstlisting}

Non-twisted affine root systems can be created as the affine extensions of finite root systems, e.g.
\begin{lstlisting}
  b2affine = makeAffineExtension[b2]
\end{lstlisting}
In the notebook interface this can be written simply as $\hat{B}_2$.

Semisimple Lie algebras can be created as the sums of simple ones:
\begin{lstlisting}[mathescape=true]
  $A_1\oplus A_1$ == finiteRootSystem[2, 2, {finiteWeight[2, {1, 0}], finiteWeight[2, {0, 1}]}]
\end{lstlisting}
%% The problem is here

The predicate \lstinline{rootSystemQ} checks if the object is a root system of finite or affine type.

The list of simple roots is a property of the root system so it is accessed as the field \lstinline{rs[simpleRoots]}.

We have implemented several functions to get major properties of root systems. The Weyl vector is computed by the function \lstinline{rho[rs_?rootSystemQ]}:
\begin{lstlisting}[label=list:1]
  In[1]  =  rho[b2]
  Out[1]  =  finiteWeight[2, {3/2, 1/2}]
\end{lstlisting}
The list of positive roots can be constructed with the function \lstinline{positiveRoots[rs_?rootSystemQ]}. For an affine Lie algebra this and related functions return the list up to some fixed grade. This grade limit is set as the value of the field \lstinline{rs[gradeLimit]} which is equal to 10 by default. The list of roots (up to \lstinline{gradeLimit}) is returned by the function \lstinline{roots[rs]}. The Cartan matrix and the fundamental weights are calculated by the functions \lstinline{cartanMatrix} and \lstinline{fundamentalWeights} correspondingly.

It is possible to specify the weight of a Lie algebra by its Dynkin labels
\begin{lstlisting}
  weight[b2][1,2] == makeFiniteWeight[{2, 1}]
\end{lstlisting}
The function \lstinline{dynkinLabels[rs_?rootSystemQ][wg_?weightQ]} returns Dynkin labels of a weight \lstinline{wg} in the root system \lstinline{rs}.

An element of the Weyl group can be specified by the set of the indices of the reflections, so the element  $w=s_{1}s_{2}s_{1}$ of the Weyl group of algebra $B_{2}$ is constructed with \lstinline{weylGroupElement[b2][1,2,1]}. Then it can be applied to the weights:
\begin{lstlisting}
  w = weylGroupElement[b2][1,2,1];
  w @ makeFiniteWeight[{1,0}] == makeFiniteWeight[{-1,0}]
\end{lstlisting}

A computation of a lexicographically minimal form \cite{casselman1994machine,casselman1995automata} for the Weyl group elements can be conveniently implemented using pattern-matching in {\it Mathematica}. In \cite{KallenShortlex} the rewrite rules for simple finite dimensional and affine Lie algebras are presented as the {\it Mathematica} patterns. Our presentation of the Weyl group elements is compatible with the code of \cite{KallenShortlex}:
\begin{lstlisting}[mathescape=true]
  In[1]  = $<<$A3reduce;
           reduce[s[1,2,1,2,1,3,2,1,1]]
  Out[1] = s[2, 3, 2]

  In[2]  = (weylGroupElement[$A_{3}$] @@ reduce[s[1,2,1,2,1,3,2,1,1]]) @ weight[$A_{3}$][-1,-2,-1]
  Out[2] = finiteWeight[4, {-2, 2, 1, -1}]

  In[3]  = dynkinLabels[$A_{3}$][Out[2]]
  Out[3] = {-4, 1, 2}
\end{lstlisting}

\subsection{Formal elements}
\label{sec:formal-elements}

We represent formal characters of modules by a special structure \lstinline{formalElement}.  The constructor \lstinline[mathescape=true]!makeFormalElement[{$\gamma_{1},\dots,\gamma_{n}$},{$m_{1},\dots,m_{n}$}]! creates a datastructure which represents the element $\sum_{i=1}^{n} m_{i} e^{\gamma_{i}}$ of the formal algebra $\mathcal{E}$. This structure is a hash-table implemented with \lstinline{DownValues}. The keys are weights at the exponents and the values are the  corresponding multiplicities. The operations in $\mathcal{E}$ are implemented for the \lstinline{formalElement} data-type: formal elements can be added, multiplied by a number or by an exponent of a weight. There exists also a multiplication of formal elements but no division.
\begin{lstlisting}[mathescape=true]
  In[1]  = makeFormalElement[{makeFiniteWeight[{1,1}],makeFiniteWeight[{0,0}]},{1,2}] *
             (2 * Exp[makeFiniteWeight[{1,0}]] *
             makeFormalElement[{makeFiniteWeight[{1,1}],makeFiniteWeight[{0,0}]},{1,2}]);
  In[2]  = In[1][weights]
  Out[2] = {finiteWeight[2, {1, 0}], finiteWeight[2, {2, 1}], finiteWeight[2, {3, 2}]}

  In[3]  = In[1][multiplicities]
  Out[3] = {8, 8, 2}
\end{lstlisting}

\subsection{Modules}
\label{sec:modules}

{\bf Affine.m} can be used to study different kinds of modules, i.e. Verma modules, irreducible modules and parabolic Verma modules.  We need the datastructure \lstinline{module} to represent a generic module of a Lie algebra $\gf$. Module properties can be deduced from its set of singular weights using the Weyl character formulae \eqref{eq:11},\eqref{eq:12},\eqref{eq:18},\eqref{eq:13}. A set of singular weights can have Weyl symmetry. It can be a symmetry with respect to the Weyl group $W_{\gf}$ or with respect to some subgroup $W_{\af}$ as in the case of parabolic Verma modules. Then it is possible to study only the main Weyl chamber $C_{\af}$. To use this symmetry a generic constructor for the \lstinline{module} datastructure accepts several parameters \lstinline{makeModule[rs_?rootSystemQ][singWeights_formalElement,subs_?rootSystemQ|emptyRootSystem[],limit:10}. Here \lstinline{rs} is the root system of a Lie algebra $\gf$, \;\lstinline{singWeights} is the set of singular weights,  \lstinline{subs} is the root system corresponding to the Weyl group $W_{\af}$ which is the (anti-)symmetry of the set of singular weights. The parameter \lstinline{limit} limits the computation for infinite-dimensional representations such as Verma or parabolic Verma modules.
There are several specialized constructors for different types of highest weight modules:
\begin{lstlisting}[mathescape=true]
vm=makeVermaModule[$B_{2}$][{2,1}];
pm=makeParabolicVermaModule[$B_{2}$][weight[$B_{2}$][2,1],{1}];
im=makeIrreducibleModule[$B_{2}$][2,1];
GraphicsRow[textPlot/@{im,vm,pm}]
$\includegraphics[width=130mm]{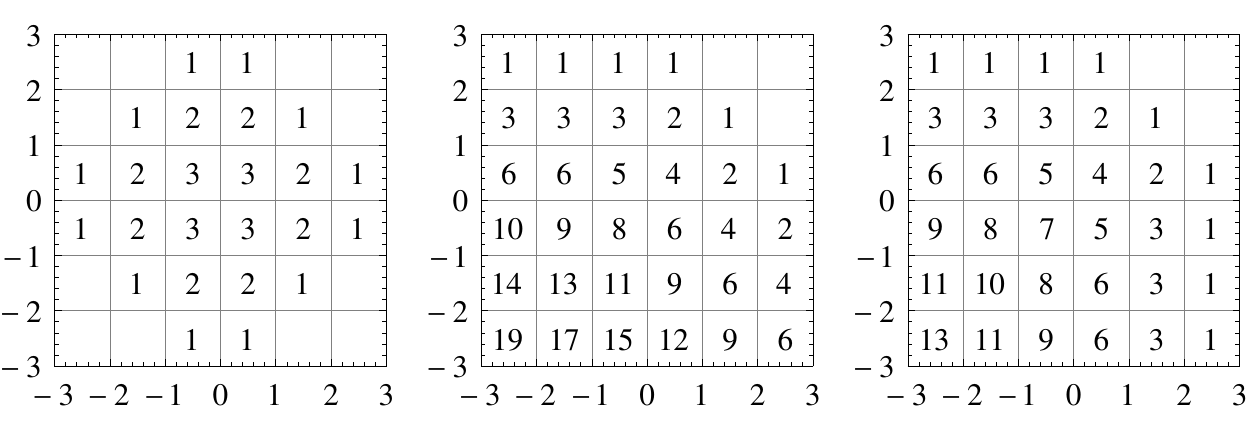}$
\end{lstlisting}

As we have already stated the properties of a module are encoded by its singular element. The function \lstinline{singularElement[m_module]} returns the singular element of a module as a \lstinline{formalElement} datastructure. The character (up to \lstinline{limit} for (parabolic) Verma modules) is returned by the function \lstinline{character[m_module]}. A direct sum of modules is a module and we use natural notation
\begin{lstlisting}[mathescape=true]
In[1] := im1=makeIrreducibleModule[$B_{2}$][weight[$B_{2}$][2,1]];
         im2=makeIrreducibleModule[$B_{2}$][weight[$B_{2}$][1,2]];
         Head[im1$\oplus$ im2]
Out[1] = module
In[2] := textPlot[im1$\oplus$ im2]
$\includegraphics[width=60mm]{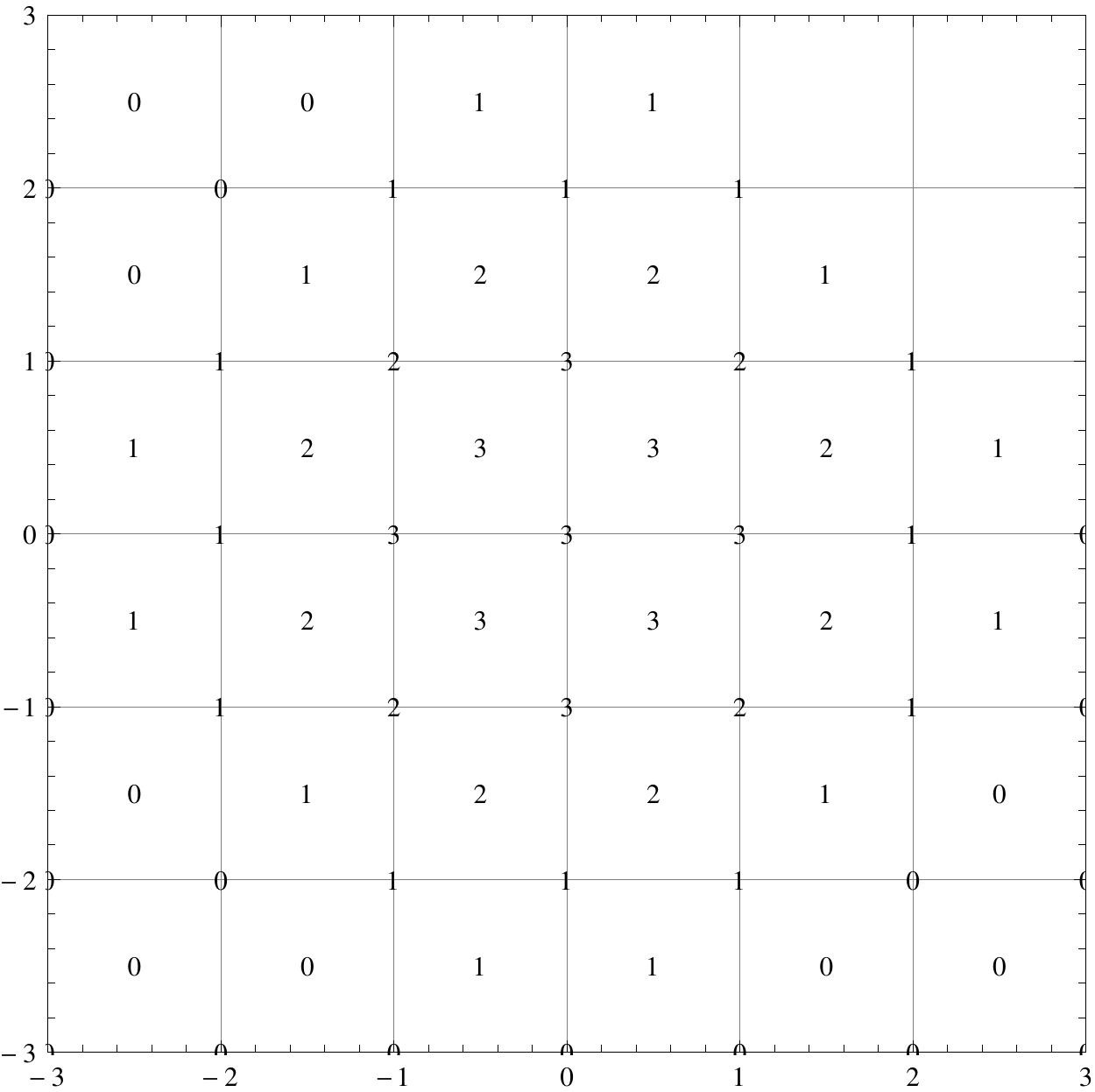}$
\end{lstlisting}

The tensor product is also implemented but only for finite-dimensional Lie algebras, since the tensor product of affine Lie algebra modules leads to rich new structures \cite{kazhdan1994tensor3,kazhdan1993tensor1,kazhdan1993tensor2} which are out of the scope of the present paper.
\begin{lstlisting}[mathescape=true]
textPlot[makeIrreducibleModule[$A_{1}$][5]$\otimes$ makeIrreducibleModule[$A_{1}$][3]];
$\includegraphics[width=60mm]{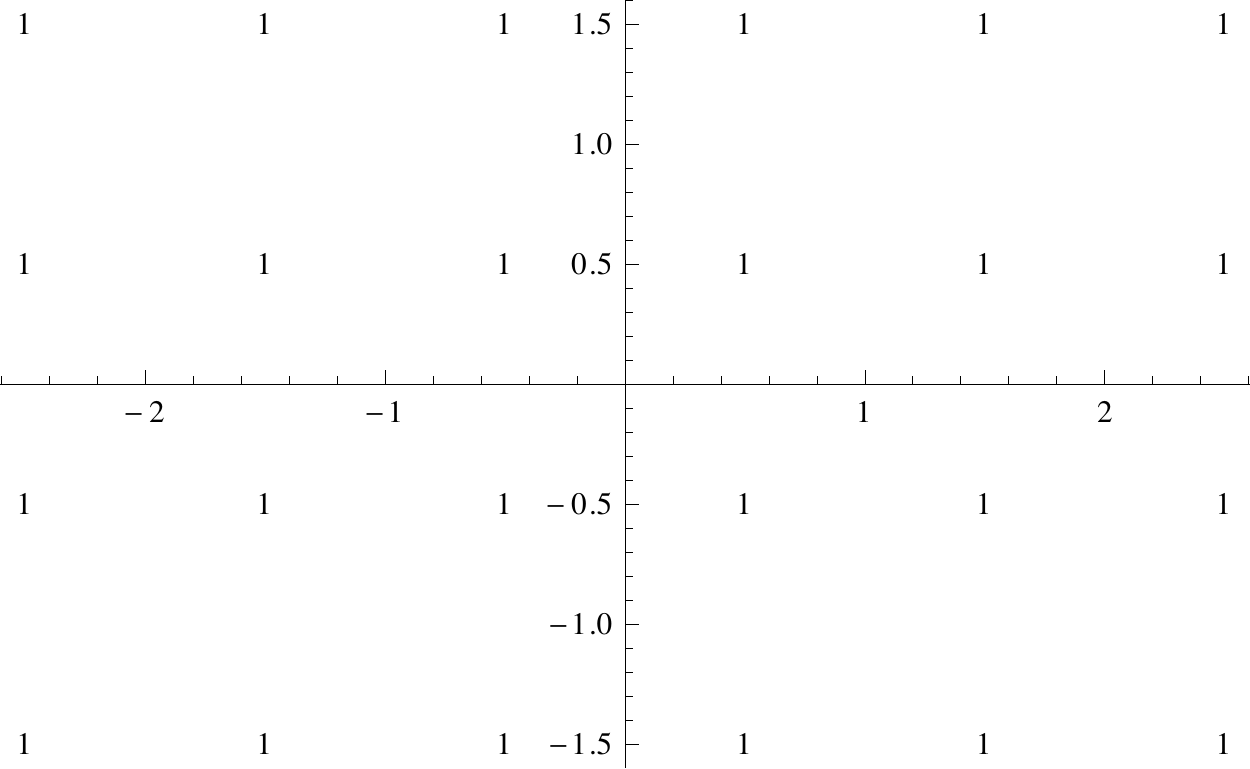}$
\end{lstlisting}
\section{Computational algorithms}
\label{sec:comp-algor}
   
As we have already stated in the section \ref{sec:high-weight-modul} there exist two recurrent relations which can be used to calculate weight multiplicities in irreducible modules. Both algorithms proceed in the following way to calculate weight multiplicities:
\begin{enumerate}
\item Create the list of weights in the main Weyl chamber by subtracting all possible combinations of simple roots from the highest weight (e.g. for a finite-dimensional algebra subtract $\alpha_{1}$ from $\mu$ while inside $\bar C$, then subtract $\alpha_{2}$ from all the weights already obtained etc).
\item Sort the list of weights by their product with the Weyl vector.
\item Use a recurrent formula. If the weight required for the recurrent computation is outside the main chamber use the Weyl symmetry.
\end{enumerate}
The difference in the performance of algorithms comes from the number of previous values required to compute the multiplicity of a weight under consideration. For  the  recurrent relation \eqref{eq:14} based on the Weyl formula it is constant and equal to the number of elements in the Weyl group (if we are far from the boundary of the representation diagram). When the Freudenthal formula \eqref{eq:15} is used the number of previous values grows with the distance from the external border of representation. So the Freudenthal formula is faster if the weight is close to the border or the rank of the algebra and the size of the Weyl group is large \cite{moody1982fast}.
Note that the Freudenthal formula is valid for the irreducible modules only, so it can not be used to study (generalized) Verma modules.

\begin{figure}[h]
  \noindent\centering{
    \includegraphics[width=80mm]{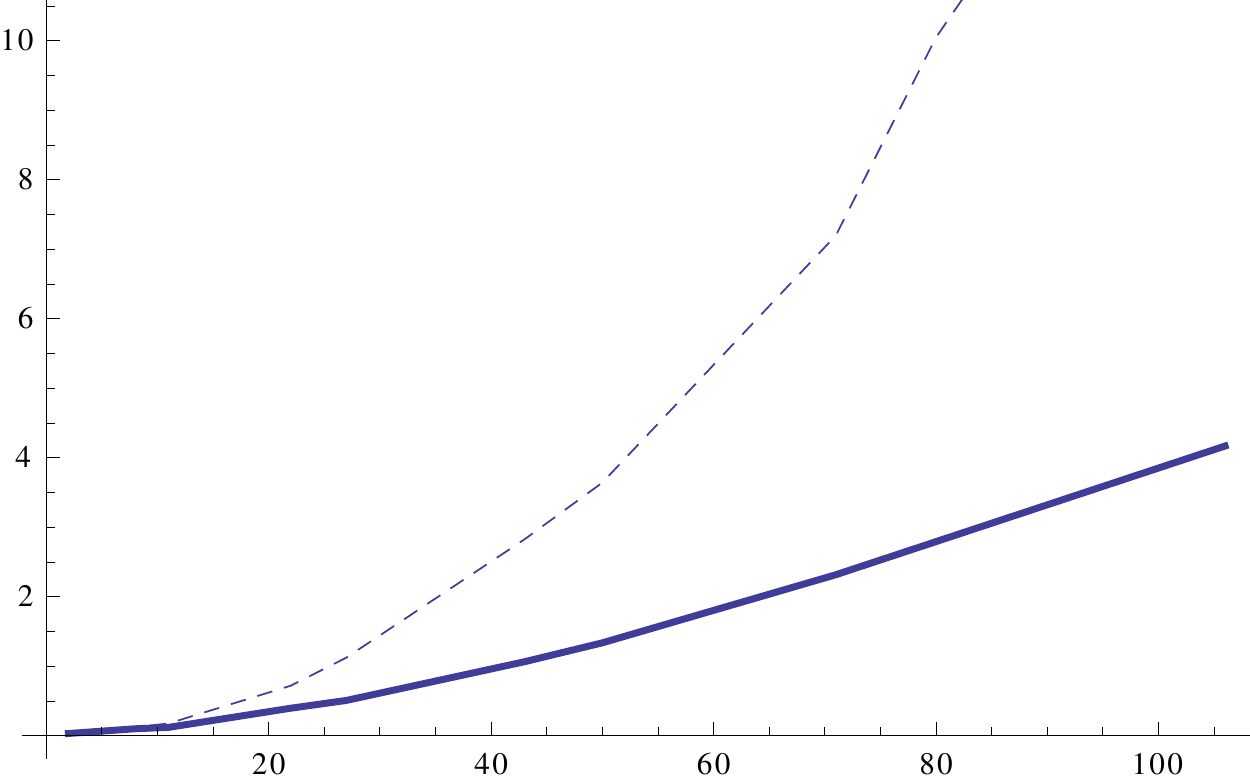}
  }
  \caption{The dependence of the running time of the algorithms based on the Freudenthal formula \eqref{eq:15} (dashed) and the recurrent relation \eqref{eq:14} (solid) on the number of weights in $\bar C$ for calculation of multiplicities in representations of $B_{2}$.}
\label{fig:freudenthal-racah-times}

\end{figure}

We have made some experiments with our implementations of the Freudenthal formula and formula \eqref{eq:15} and prepared the Fig.  \ref{fig:freudenthal-racah-times}, which depicts the dependence of the computation time on the number of weights in a module.

In the calculation of branching coefficients the application of the Freudenthal formula requires a complete construction of the formal characters of an algebra module and all the representations of a subalgebra. It is impractical if the ranks of an algebra and a subalgebra are big, for example for the maximal subalgebras.

An alternative algorithm was presented in the paper \cite{2010arXiv1007.0318L}. It contains the following steps:

\begin{enumerate}
\item  Construct the root system $\Delta _{\af}$ for an embedding $\af\rightarrow \gf$.
\item  Select all the positive roots $\alpha \in \Delta ^{+}$ orthogonal to  $\af$, i.e. form the set $\Delta_{\afb }^{+}$.
\item  Construct the set $\Gamma _{\af\rightarrow \gf}$. The relation (\ref{eq:6}) defines the sign function
 $s(\gamma)$ and the set $\Phi_{\af\subset \frak{g}}$. The lowest weight
 $\gamma_0$ is subtracted to get the fan:
 $\Gamma _{\af\rightarrow \frak{g}}=\left\{ \xi -\gamma _{0}|\xi \in \Phi _{%
\af\subset \frak{g}}\right\} \setminus \left\{ 0\right\}$.

\item  Construct the set $\widehat{\Psi ^{(\mu )}}=\left\{ w (\mu +\rho
)-\rho ;\;w \in W\right\} $ of singular weights for the $\frak{g}$%
-module $L^{(\mu )}$.

\item  Select the weights $\left\{ \mu _{\widetilde{\af_{\perp }}%
}\left( w\right) =\pi _{\widetilde{\af_{\perp }}}\left[ w(\mu +\rho
)-\rho \right] -\mathcal{D}_{\af_{\perp }}\in \overline{C_{\widetilde{%
\af_{\perp }}}}\right\} $. Since the set $\Delta_{\afb }^{+}$ is fixed
we can easily check that the weight $\mu _{\widetilde{\af_{\perp }}%
}\left( w\right) $ belongs to the main Weyl chamber $\overline{C_{\widetilde{%
\af_{\perp }}}}$ (by computing its scalar product with the fundamental
weights of $\afb$).

\item  For the weights $\mu _{\widetilde{\af_{\perp }}}\left( w\right) $
calculate dimensions of the corresponding modules, $\mathrm{\dim }\left(
L_{\widetilde{\af_{\perp }}}^{\mu _{\widetilde{\af_{\perp }}%
}\left( u\right) }\right) $, using the Weyl dimension formula and construct
the singular element $\Psi ^{\left( \mu \right) }_{\left(  \af, \afb \right)}$.

\item  Calculate the signed branching coefficients using 
the recurrent relation (\ref{recurrent-rel}) and select among them those
corresponding to the weights in the main Weyl
chamber $\overline{C_{\af}}$.
\end{enumerate}

We can speed up the algorithm by
one-time computation of the representatives of conjugate classes $W/W_{\afb }$.

Consider the regular embedding $B_{2}\subset B_{4}$. In this case the fan consists of 24 elements. In order to decompose $B_{4}$ module we need to construct the subset of singular weights of the module which projects to the main Weyl chamber of the subalgebra $B_{2}$. The full set of singular weights consists of 384 elements. The required subset contains at most 48 elements. The time for the construction of this required subset is negligible if the number of branching coefficients is greater than that.

 We may estimate the total number of required operations for the computation of branching coefficients as the product of the number of elements with non-zero branching coefficients  in the main Weyl chamber of a subalgebra  and the number of elements in the fan. So we have a linear growth.  If we use a direct algorithm we need to compute the multiplicities for each module in the decomposition. So the number of operations grows faster than the square of the number of elements with non-zero branching coefficients in the main Weyl chamber of a subalgebra. 

\begin{figure}[h]
  \noindent\centering{
   \includegraphics[width=100mm]{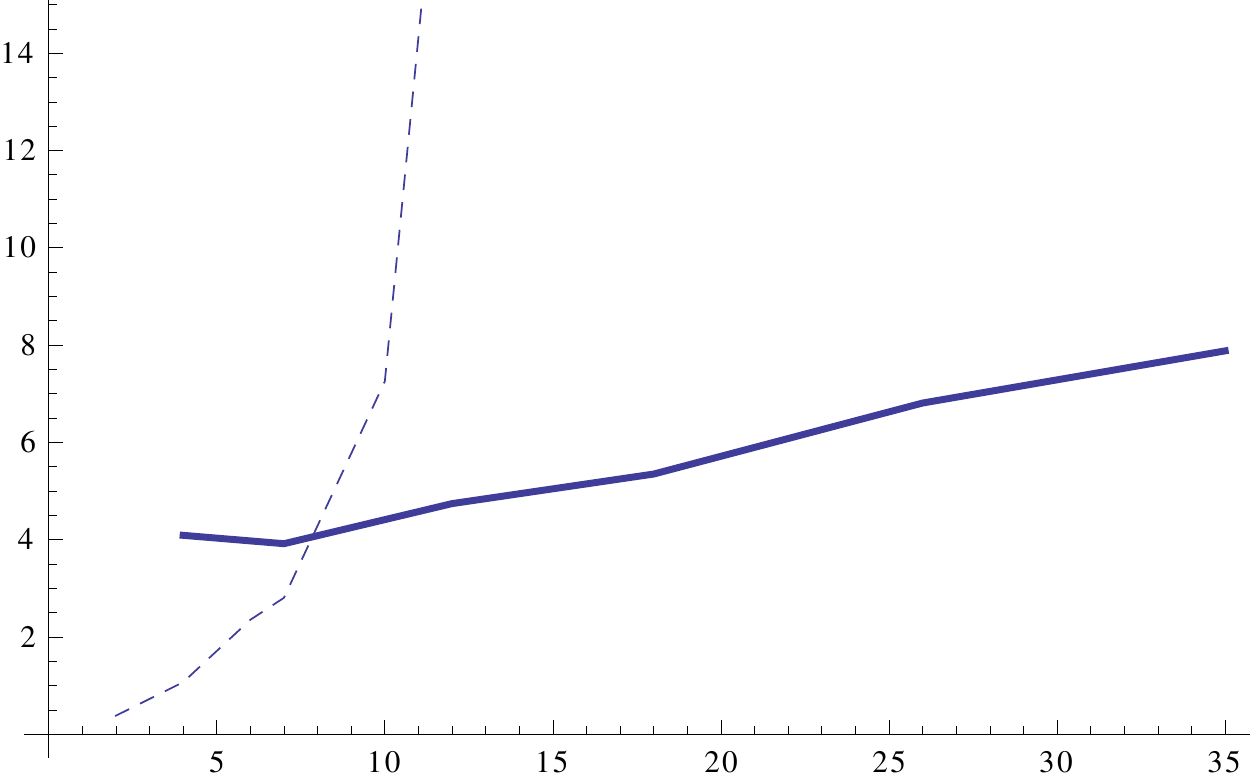}
  }
  \caption{Dependence of the time required to compute branching coefficients for $B_{3}\subset B_{4}$ on  the number of weights in $\bar C$. The dashed line corresponds to the direct  algorithm based on the Freudenthal formula \eqref{eq:15}, the solid one to the use of recurrent relation \eqref{recurrent-rel}.}
  \label{fig:branching}
\end{figure}

To illustrate this performance issue we present the Figure \ref{fig:branching} where we show the time required to compute the branching coefficients for $B_{3}\subset B_{4}$.

\section{Examples}
\label{sec:examples}
In this section we present some examples of computations available with {\bf Affine.m} and the code required to produce these results.

\subsection{Tensor product decompositon for finite-dimensional Lie algebras}
\label{sec:tens-prod-decomp}

Computation of fusion coefficients for decomposition of the tensor product of highest-weight modules to the direct sum of irreducible modules has numerous applications in physics. For example, we can consider the spin of a composite system such as an atom. Another interesting example is the integrable spin chain consisting of $N$ particles with the spins living in some representation $L$ of a Lie algebra $\gf$ with a $\gf$-invariant Hamiltonian $H$, describing nearest-neighbour spin-spin interaction. In order to solve such a system, i.e. to find eigenstates of the Hamiltonian, we need to decompose $L^{\otimes N}$ into the direct sum of the irreducible $\gf$-modules of lower dimension and diagonalize the Hamiltonian on these modules.

For the fundamental representations of simple Lie algebras it is sometimes possible to get an analytic formula for the dependence of the decomposition coefficients on $N$ (See \cite{LyakhovskyPostnova2011}). Our code provides the numerical values and can be used to check the analytic results.

Consider for example a fourth tensor power of the first fundamental representation $\left(L^{[1,0]}\right)^{\otimes 4}$ of algebra $B_{2}$ . Decomposition coefficients are just the branching
 coefficients for the tensor power module reduced to the diagonal subalgebra $B_{2}\subset B_{2}\oplus B_{2}\oplus B_{2}\oplus B_{2}$. So the following code calculates these coefficients:
\begin{lstlisting}[mathescape=true]
fm = makeIrreducibleModule[$B_{2}$][1, 0];
tp = ((fm$\otimes$ ]fm)$\otimes$ fm)$\otimes$]fm;
subs = makeFiniteRootSystem[
  {1/4*{1, -1, 1, -1, 1, -1, 1, -1}, 
   1/4*{0, 1, 0, 1, 0, 1, 0, 1}}];
bc = branching[tp, subs];
{bc[#], dynkinLabels[subs][#]} & /@ bc[weights]
\end{lstlisting}
It produces a list of highest weights and tensor product decomposition coefficients:
\begin{lstlisting}
{{1, {4, 0}}, {3, {2, 2}}, {0, {3, 0}}, 
{2, {0, 4}}, {3, {1, 2}}, {6, {2, 0}}, 
{6, {0, 2}}, {1, {1, 0}}, {3, {0, 0}}}]
\end{lstlisting}

Returning to the problem of spin chain Hamiltonian diagonalization we can see that instead of diagonalizing the operator in a space of dimension $625$ we can diagonalize the operators in the spaces of dimensions $55, 81, 30, 35, 35, 14, 10, 5, 1$.

\subsection{Branching and parabolic Verma modules}
\label{sec:branch-parab-verma}

\begin{figure}[h]
  \noindent\centering{
    \includegraphics[width=80mm]{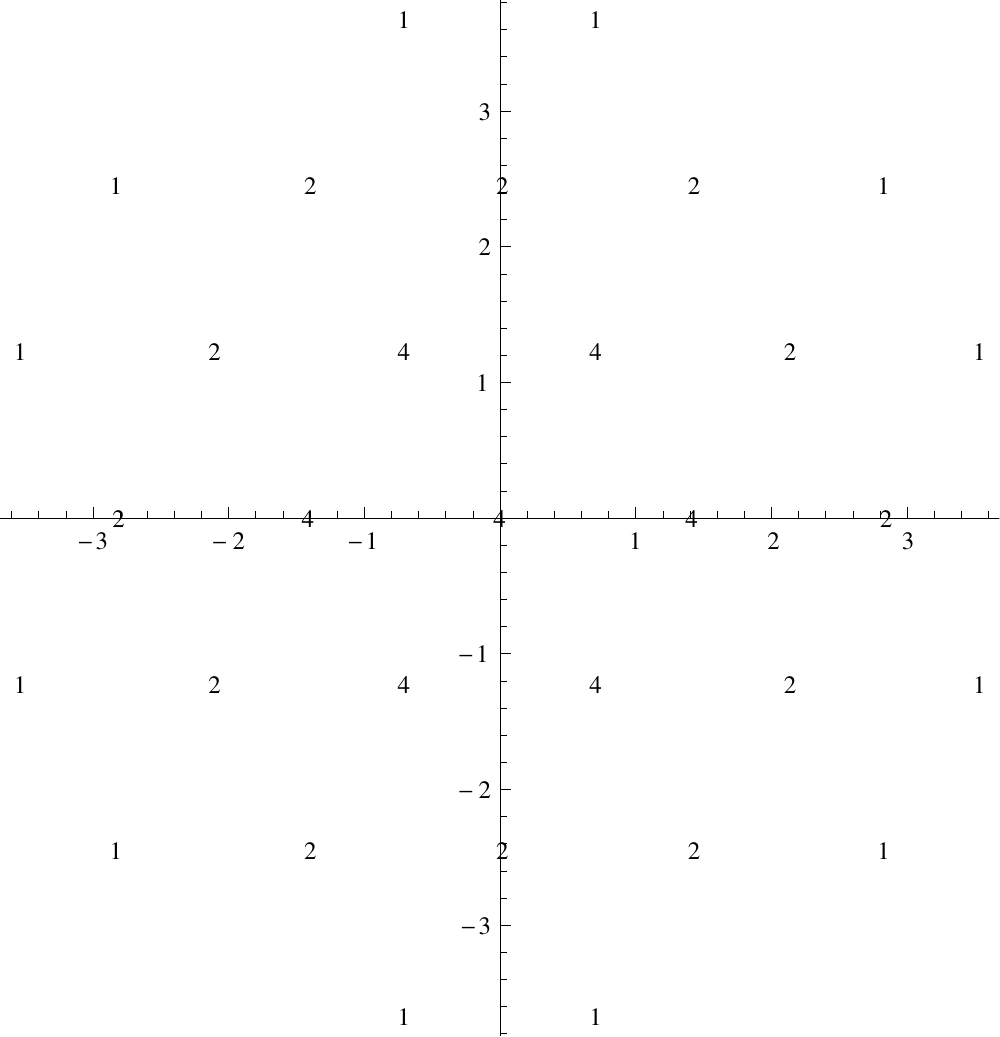}
  }
  \caption{Character of the irreducible $G_{2}$-module $L^{[1,1]}$}
  \label{branching-bgg}
\end{figure}

\begin{figure}[h]
  \noindent\centering{
    \includegraphics[width=150mm]{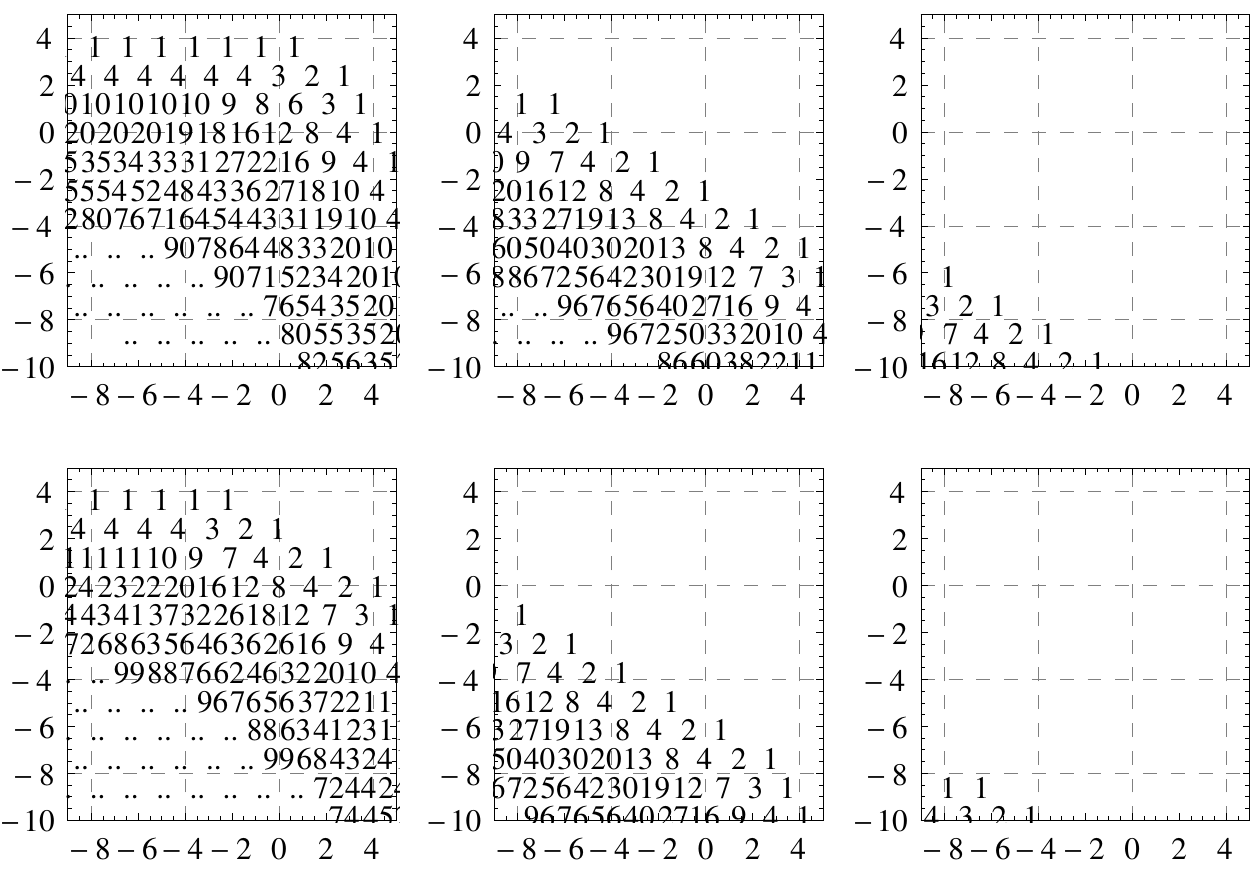}
  }
  \caption{The characters of  generalized Verma modules of $G_{2}$ appearing in the decomposition of $L^{[1,1]}$. The parabolic Verma modules in the upper row appear in the decomposition with a positive sign, in  lower row with a negative.}
  \label{g2-pverma}
\end{figure}

We illustrate the generalized BGG-resolution by the diagrams of $G_{2}$ parabolic Verma modules which appear in the decomposition of the irreducible module $L^{[1,1]}_{G_{2}}$:
\begin{equation}
\mathrm{ch}\left( L^{\mu }\right) =\sum_{u\in U}\;e^{\mu _{\aft}\left(
u\right) }\epsilon (u)\mathrm{ch}M_{I}^{\mu _{\frak{a}_{\perp }}\left(
u\right) }.  \label{char-in-gen-verma-mod}
\end{equation}
The character of $L^{[1,1]}$ is presented in Figure \ref{branching-bgg}, the characters of the generalized Verma modules in the decomposition \eqref{char-in-gen-verma-mod} are shown in Figure \ref{g2-pverma}. The characters in the upper row appear in \eqref{char-in-gen-verma-mod} with a positive sign and in the lower row with a negative.


\subsection{String functions of affine Lie algebras and CFT models}
\label{sec:string-funct-affine}

String functions can be used to present a formal character of an affine Lie algebra highest weight representation. They have interesting analytic and modular properties \cite{kac1990idl,kac1988modular,kac1984infinite}.

{\bf Affine.m} produces power series decomposition for string functions. Consider an affine Lie algebra $\hat{sl(3)}=\hat A_{2}$ and its highest weight module $L^{(1,0,0)}$. To get the string functions we can use the code:
\begin{lstlisting}[mathescape=true]
stringFunctions[$\hat A_2$,{1,1,2}]
{{0, 0, 4}, 
  $2 q + 10 q^2 + 40 q^3 + 133 q^4 + 398 q^5 + 1084 q^6 + 2760 q^7 + 6632 q^8 + 15214 q^9 + 33508 q^{10}$}, 
{{0, 3, 1}, 
  $2 q + 12 q^2 + 49 q^3 + 166 q^4 + 494 q^5 + 1340 q^6 + 3387 q^7 + 8086 q^8 + 18415 q^9 + 40302 q^{10}$}, 
{{1, 1, 2}, 
  $1 + 6 q + 27 q^2 + 96 q^3 + 298 q^4 + 836 q^5 + 2173 q^6 + 5310 q^7 + 12341 q^8 + 27486 q^9 + 59029 q^{10}$}, 
{{2, 2, 0}, 
  $1 + 8 q + 35 q^2 + 124 q^3 + 379 q^4 + 1052 q^5 + 2700 q^6 + 6536 q^7 + 15047 q^8 + 33248 q^9 + 70877 q^{10}$}, 
{{3, 0, 1}, 
  $2 + 12 q + 49 q^2 + 166 q^3 + 494 q^4 + 1340 q^5 + 3387 q^6 + 8086 q^7 + 18415 q^8 + 40302 q^9 + 85226 q^{10}$}
\end{lstlisting}
Similarly for the affine Lie algebra $\hat G_{2}$ we get
\begin{lstlisting}[mathescape=true]
stringFunctions[$\hat G_2$,{1,1,0}]
{{2, 0, 0}, 
  $1 + 8 q + 37 q^2 + 138 q^3 + 431 q^4 + 1227 q^5 + 3208 q^6 + 7901 q^7$}, 
{{0, 0, 1}, 
  $3 q + 18 q^2 + 73 q^3 + 247 q^4 + 736 q^5 + 2000 q^6 + 5070 q^7$}, 
{{1, 1, 0},
  $1 + 7 q + 32 q^2 + 117 q^3 + 370 q^4 + 1055 q^5 + 2780 q^6 + 6880 q^7$}, 
{{0, 2, 0}, 
  $3 q + 15 q^2 + 63 q^3 + 210 q^4 + 633 q^5 + 1725 q^6 + 4407 q^7$}
\end{lstlisting}

\subsection{Branching functions and coset models of conformal field theory}
\label{sec:branch-funct-coset}

It is believed that most rational models of CFT can be obtained from the cosets $G/A$ corresponding to the embedding $\af\subset\gf$. These models can be studied as gauge theories \cite{Hwang:1994yr, hwang1993brst}.

Branching functions for an embedding $\af\subset\gf$ are the partition functions of CFT on the torus (see \cite{difrancesco1997cft}).

As a first example we show how to construct branching functions for the embedding $\hat A_{1}\to \hat B_{2}$ up to the tenth grade:
\begin{lstlisting}[mathescape=true]
branchingFunctions[$\hat B_{2}$,makeAffineExtension[makeFiniteRootSystem[{{1, 1}}]], {1, 1, 1}]

 {{3, 0}, 
  $2 + 14 q + 52 q^2 + 154 q^3 + 410 q^4 + 994 q^5 + 2248 q^6 + 4832 q^7 + 9934 q^8 + 19680 q^9 + 37802 q^{10}$},
 {{2, 1}, 
  $4 + 20 q + 72 q^2 + 220 q^3 + 584 q^4 + 1424 q^5 + 3248 q^6 + 7012 q^7 + 14488 q^8 + 28844 q^9 + 55616 q^{10}$},
 {{0, 3}, 
  $4 q + 20 q^2 + 68 q^3 + 200 q^4 + 516 q^5 + 1224 q^6 + 2736 q^7 + 5808 q^8 + 11820 q^9 + 23236 q^{10}$},
 {{1, 2}, 
  $2 + 14 q + 54 q^2 + 168 q^3 + 462 q^4 + 1148 q^5 + 2656 q^6 + 5812 q^7 + 12130 q^8 + 24358 q^9 + 47328 q^{10}$}
\end{lstlisting}

Another example demonstrates the computation of branching functions for the regular embedding $\hat B_{2}\subset \hat C_{3}$:
\begin{lstlisting}[mathescape=true]
sub=makeAffineExtension[parabolicSubalgebra[$C_{3}$][2,3]];
branchingFunctions[$\hat C_{3}$,sub, {2, 0, 0, 0}]

{{0, 1, 0}, 
  $2 q - 20 q^3 + 24 q^4 + 82 q^5 - 320 q^6 + 108 q^7$}, 
{{1, 0, 0}, 
  $1 - q - 8 q^2 + 19 q^3 + 16 q^4 - 156 q^5 + 205 q^6 + 640 q^7$}, 
{{0, 0, 1}, 
  $q - 5 q^3 + 7 q^5$}
\end{lstlisting}

\section{Conclusion}
\label{sec:conclusion}
We have presented the package {\bf Affine.m} for computations in representation theory of finite-dimensional and affine Lie algebras. It can be used to study Weyl symmetry, root systems, irreducible, Verma and parabolic Verma modules of finite-dimensional and affine Lie algebras.  In the present paper we have also discussed main ideas used for the implementation of the package and described the most important notions of representation theory required to use {\bf Affine.m}. 

We have demonstrated that the recurrent approach based on the Weyl character formula is not only useful for calculations but also allows us to establish connections with the (generalized) Bernstein-Gelfand-Gelfand resolution. 

Also we have presented examples of computations connected with problems of physics and mathematics. 

In future versions of our software we are going to treat twisted affine Lie algebras, extended affine Lie algebras and provide more direct support for tensor product decompositions. 

\section*{Acknowledgements}
\label{sec:acknowledgements}
I thank V.D. Lyakhovsky and O.Postnova for discussions and helpful comments.

The work is supported by the Chebyshev Laboratory
(Department of Mathematics and Mechanics, Saint-Petersburg State
University) under the grant 11.G34.31.0026 of the Government of the
Russian Federation.

\appendix

\section{Software package}
\label{package}
The package can be freely downloaded from \url{http://github.com/naa/Affine}. To get the development code use the command
\begin{lstlisting}[language=bash]
 git clone git://github.com/naa/Affine.git
\end{lstlisting}

Contents of the package:
\begin{verbatim}
    Affine/                                root folder
      demo/                                  demonstrations
        demo.nb                                demo notebook
        paper.nb                               code for the paper
      doc/                                 documentation folder
        figures/                             figures in paper 
          timing.pdf                           diagram showing performance
          branching-timing.pdf                 ...  for branching coefficients  
          irrep-sum.pdf                        sum of B2 irreps
          irrep-verma-pverma.pdf               irrep, Verma, (p)Verma for B2
          G2-irrep.pdf                         irrep for G2
          G2-pverma.pdf                        parabolic Verma for G2
          tensor-product.pdf                   tensor product of A1-modules
        bibliography.bib                     bibliographic database
        paper.pdf                            present paper
        paper.tex                            paper source
        TODO.org                             list of issues
      src/                                 source folder
        affine.m                             main software package
      tests/                               unit tests folder
        tests.m                              unit tests
      README.markdown                      installation and usage notes
\end{verbatim}

%% References


\section*{References}
\label{sec:references}

\bibliography{bibliography}

\begin{thebibliography}{44}
\expandafter\ifx\csname natexlab\endcsname\relax\def\natexlab#1{#1}\fi
\providecommand{\bibinfo}[2]{#2}
\ifx\xfnm\relax \def\xfnm[#1]{\unskip,\space#1}\fi
%Type = Book
\bibitem[{Belinfante and Kolman(1989)}]{belinfante1989survey}
\bibinfo{author}{J.~Belinfante}, \bibinfo{author}{B.~Kolman},
  \bibinfo{title}{{A survey of Lie groups and Lie algebras with applications
  and computational methods}}, \bibinfo{publisher}{Society for Industrial
  Mathematics}, \bibinfo{year}{1989}.
%Type = Article
\bibitem[{Nutma(2011)}]{simplie}
\bibinfo{author}{T.~Nutma},
\newblock \bibinfo{title}{Simplie}  (\bibinfo{year}{2011}),
  \url{http://code.google.com/p/simplie/}.
%Type = Article
\bibitem[{Van~Leeuwen(1994)}]{vanleeuwen1994lsp}
\bibinfo{author}{M.~Van~Leeuwen},
\newblock \bibinfo{title}{{LiE, a software package for Lie group
  computations}},
\newblock \bibinfo{journal}{Euromath Bull} \bibinfo{volume}{1}
  (\bibinfo{year}{1994}) \bibinfo{pages}{83--94},
  \url{http://www-math.univ-poitiers.fr/~maavl/LiE/}.
%Type = Article
\bibitem[{Stembridge(1995)}]{stembridge1995mps}
\bibinfo{author}{J.~Stembridge},
\newblock \bibinfo{title}{{A Maple package for symmetric functions}},
\newblock \bibinfo{journal}{Journal of Symbolic Computation}
  \bibinfo{volume}{20} (\bibinfo{year}{1995}) \bibinfo{pages}{755--758}.
%Type = Article
\bibitem[{Stembridge(2011)}]{coxweyl}
\bibinfo{author}{J.~Stembridge},
\newblock \bibinfo{title}{Coxeter/weyl packages for maple}
  (\bibinfo{year}{2011}),
  \url{{http://www.math.lsa.umich.edu/~jrs/maple.html}}.
%Type = Article
\bibitem[{Fischbacher(2002)}]{fischbacher2002ilp}
\bibinfo{author}{T.~Fischbacher},
\newblock \bibinfo{title}{{Introducing LambdaTensor1.0 - A package for explicit
  symbolic and numeric Lie algebra and Lie group calculations}}
  (\bibinfo{year}{2002}), \href{http://arxiv.org/abs/hep-th/0208218}{{\tt
  arXiv:hep-th/0208218}}.
%Type = Article
\bibitem[{Fuchs et~al.(1996)Fuchs, Schellekens, and Schweigert}]{Fuchs:1996dd}
\bibinfo{author}{J.~Fuchs}, \bibinfo{author}{A.~N. Schellekens},
  \bibinfo{author}{C.~Schweigert},
\newblock \bibinfo{title}{{A matrix S for all simple current extensions}},
\newblock \bibinfo{journal}{Nucl. Phys.} \bibinfo{volume}{B473}
  (\bibinfo{year}{1996}) \bibinfo{pages}{323--366},
  \href{http://arxiv.org/abs/hep-th/9601078}{{\tt arXiv:hep-th/9601078}},
  \url{http://www.nikhef.nl/~t58/kac.html}.
%Type = Article
\bibitem[{Moody and Patera(1982)}]{moody1982fast}
\bibinfo{author}{R.~Moody}, \bibinfo{author}{J.~Patera},
\newblock \bibinfo{title}{Fast recursion formula for weight multiplicities},
\newblock \bibinfo{journal}{Bulletin (New Series) of the American Mathematical
  Society} \bibinfo{volume}{7} (\bibinfo{year}{1982})
  \bibinfo{pages}{237--242}.
%Type = Article
\bibitem[{Stembridge(2001)}]{stembridge2001computational}
\bibinfo{author}{J.~Stembridge},
\newblock \bibinfo{title}{{Computational aspects of root systems, Coxeter
  groups, and Weyl characters, Interaction of combinatorics and representation
  theory, MSJ Mem., vol. 11}},
\newblock \bibinfo{journal}{Math. Soc. Japan, Tokyo}  (\bibinfo{year}{2001})
  \bibinfo{pages}{1--38}.
%Type = Article
\bibitem[{Casselman(1994)}]{casselman1994machine}
\bibinfo{author}{W.~Casselman},
\newblock \bibinfo{title}{{Machine calculations in Weyl groups}},
\newblock \bibinfo{journal}{Inventiones mathematicae} \bibinfo{volume}{116}
  (\bibinfo{year}{1994}) \bibinfo{pages}{95--108}.
%Type = Book
\bibitem[{Kac(1990)}]{kac1990idl}
\bibinfo{author}{V.~Kac}, \bibinfo{title}{{Infinite dimensional Lie algebras}},
  \bibinfo{publisher}{Cambridge University Press}, \bibinfo{year}{1990}.
%Type = Article
\bibitem[{Walton(1999)}]{Walton:1999xc}
\bibinfo{author}{M.~Walton},
\newblock \bibinfo{title}{{Affine Kac-Moody algebras and the Wess-Zumino-Witten
  model}}  (\bibinfo{year}{1999}),
  \href{http://arxiv.org/abs/hep-th/9911187}{{\tt arXiv:hep-th/9911187}}.
%Type = Book
\bibitem[{Di~Francesco et~al.(1997)Di~Francesco, Mathieu, and
  Senechal}]{difrancesco1997cft}
\bibinfo{author}{P.~Di~Francesco}, \bibinfo{author}{P.~Mathieu},
  \bibinfo{author}{D.~Senechal}, \bibinfo{title}{{Conformal field theory}},
  \bibinfo{publisher}{Springer}, \bibinfo{year}{1997}.
%Type = Article
\bibitem[{Goddard et~al.(1985)Goddard, Kent, and Olive}]{Goddard198588}
\bibinfo{author}{P.~Goddard}, \bibinfo{author}{A.~Kent},
  \bibinfo{author}{D.~Olive},
\newblock \bibinfo{title}{Virasoro algebras and coset space models},
\newblock \bibinfo{journal}{Physics Letters B} \bibinfo{volume}{152}
  (\bibinfo{year}{1985}) \bibinfo{pages}{88 -- 92}.
%Type = Article
\bibitem[{Dunbar and Joshi(1993)}]{Dunbar:1992gh}
\bibinfo{author}{D.~C. Dunbar}, \bibinfo{author}{K.~G. Joshi},
\newblock \bibinfo{title}{{Characters for coset conformal field theories}},
\newblock \bibinfo{journal}{Int. J. Mod. Phys.} \bibinfo{volume}{A8}
  (\bibinfo{year}{1993}) \bibinfo{pages}{4103--4122},
  \href{http://arxiv.org/abs/hep-th/9210122}{{\tt arXiv:hep-th/9210122}}.
%Type = Article
\bibitem[{Gannon(2001)}]{gannon2001algorithms}
\bibinfo{author}{T.~Gannon},
\newblock \bibinfo{title}{{Algorithms for affine Kac-Moody algebras}}
  (\bibinfo{year}{2001}), \href{http://arxiv.org/abs/hep-th/0106123}{{\tt
  arXiv:hep-th/0106123}}.
%Type = Book
\bibitem[{Kass et~al.(1990)Kass, Moody, Patera, and Slansky}]{kass1990ala}
\bibinfo{author}{S.~Kass}, \bibinfo{author}{R.~Moody},
  \bibinfo{author}{J.~Patera}, \bibinfo{author}{R.~Slansky},
  \bibinfo{title}{{Affine Lie algebras, weight multiplicities, and branching
  rules}}, \bibinfo{publisher}{Sl}, \bibinfo{year}{1990}.
%Type = Book
\bibitem[{Humphreys(1997)}]{humphreys1997introduction}
\bibinfo{author}{J.~Humphreys}, \bibinfo{title}{{Introduction to Lie algebras
  and representation theory}}, \bibinfo{publisher}{Springer},
  \bibinfo{year}{1997}.
%Type = Book
\bibitem[{Humphreys(1992)}]{humphreys1992reflection}
\bibinfo{author}{J.~Humphreys}, \bibinfo{title}{{Reflection groups and Coxeter
  groups}}, \bibinfo{publisher}{Cambridge Univ Pr}, \bibinfo{year}{1992}.
%Type = Book
\bibitem[{Wakimoto(2001{\natexlab{a}})}]{wakimoto2001idl}
\bibinfo{author}{M.~Wakimoto}, \bibinfo{title}{{Infinite-dimensional Lie
  algebras}}, \bibinfo{publisher}{American Mathematical Society},
  \bibinfo{year}{2001}{\natexlab{a}}.
%Type = Book
\bibitem[{Wakimoto(2001{\natexlab{b}})}]{wakimoto2001lectures}
\bibinfo{author}{M.~Wakimoto}, \bibinfo{title}{{Lectures on
  infinite-dimensional Lie algebra}}, \bibinfo{publisher}{World scientific},
  \bibinfo{year}{2001}{\natexlab{b}}.
%Type = Book
\bibitem[{Fulton and Harris(1991)}]{fulton1991representation}
\bibinfo{author}{W.~Fulton}, \bibinfo{author}{J.~Harris},
  \bibinfo{title}{{Representation theory: a first course}}, volume
  \bibinfo{volume}{129}, \bibinfo{publisher}{Springer Verlag},
  \bibinfo{year}{1991}.
%Type = Book
\bibitem[{Bourbaki(2002)}]{bourbaki2002lie}
\bibinfo{author}{N.~Bourbaki}, \bibinfo{title}{{Lie groups and Lie algebras}},
  \bibinfo{publisher}{Springer Verlag}, \bibinfo{year}{2002}.
%Type = Book
\bibitem[{Humphreys(2008)}]{humphreys2008representations}
\bibinfo{author}{J.~Humphreys}, \bibinfo{title}{{Representations of semisimple
  Lie algebras in the BGG category O}}, \bibinfo{publisher}{Amer Mathematical
  Society}, \bibinfo{year}{2008}.
%Type = Book
\bibitem[{Carter(2005)}]{carter2005lie}
\bibinfo{author}{R.~Carter}, \bibinfo{title}{{Lie algebras of finite and affine
  type}}, \bibinfo{publisher}{Cambridge University Press},
  \bibinfo{year}{2005}.
%Type = Article
\bibitem[{Bernstein et~al.(1976)Bernstein, Gel'fand, and
  Gel'fand}]{bernstein1976category}
\bibinfo{author}{J.~Bernstein}, \bibinfo{author}{I.~Gel'fand},
  \bibinfo{author}{S.~Gel'fand},
\newblock \bibinfo{title}{{Category of g-modules}},
\newblock \bibinfo{journal}{Funktsional'nyi Analiz i ego prilozheniya}
  \bibinfo{volume}{10} (\bibinfo{year}{1976}) \bibinfo{pages}{1--8}.
%Type = Article
\bibitem[{Bernstein et~al.(1971)Bernstein, Gel'fand, and
  Gel'fand}]{bernstein1971structure}
\bibinfo{author}{I.~Bernstein}, \bibinfo{author}{I.~Gel'fand},
  \bibinfo{author}{S.~Gel'fand},
\newblock \bibinfo{title}{{Structure of representations generated by vectors of
  highest weight}},
\newblock \bibinfo{journal}{Functional Analysis and Its Applications}
  \bibinfo{volume}{5} (\bibinfo{year}{1971}) \bibinfo{pages}{1--8}.
%Type = Article
\bibitem[{Il'in et~al.(2010)Il'in, Kulish, and Lyakhovsky}]{il2010folded}
\bibinfo{author}{M.~Il'in}, \bibinfo{author}{P.~Kulish},
  \bibinfo{author}{V.~Lyakhovsky},
\newblock \bibinfo{title}{{Folded fans and string functions}},
\newblock \bibinfo{journal}{Zapiski Nauchnykh Seminarov POMI}
  \bibinfo{volume}{374} (\bibinfo{year}{2010}) \bibinfo{pages}{197--212}.
%Type = Article
\bibitem[{Kulish and Lyakhovsky(2008)}]{kulish4sfa}
\bibinfo{author}{P.~Kulish}, \bibinfo{author}{V.~Lyakhovsky},
\newblock \bibinfo{title}{{String Functions for Affine Lie Algebras Integrable
  Modules}},
\newblock \bibinfo{journal}{Symmetry, Integrability and Geometry: Methods and
  Applications} \bibinfo{volume}{4} (\bibinfo{year}{2008}),
  \href{http://arxiv.org/abs/0812.2381}{{\tt arXiv:0812.2381 [math.RT]}}.
%Type = Article
\bibitem[{Lyakhovsky and Nazarov(2011)}]{2010arXiv1007.0318L}
\bibinfo{author}{V.~Lyakhovsky}, \bibinfo{author}{A.~Nazarov},
\newblock \bibinfo{title}{Recursive algorithm and branching for nonmaximal
  embeddings},
\newblock \bibinfo{journal}{Journal of Physics A: Mathematical and Theoretical}
  \bibinfo{volume}{44} (\bibinfo{year}{2011}) \bibinfo{pages}{075205},
  \href{http://arxiv.org/abs/1007.0318}{{\tt arXiv:1007.0318 [math.RT]}}.
%Type = Article
\bibitem[{{Lyakhovsky} and {Nazarov}(2011)}]{2011arXiv1102.1702L}
\bibinfo{author}{V.~{Lyakhovsky}}, \bibinfo{author}{A.~{Nazarov}},
\newblock \bibinfo{title}{{Recursive properties of branching and BGG
  resolution}},
\newblock \bibinfo{journal}{{Theor. Math. Phys.}} \bibinfo{volume}{{169}}
  (\bibinfo{year}{{2011}}) \bibinfo{pages}{{1551--1560}},
  \href{http://arxiv.org/abs/1102.1702}{{\tt arXiv:1102.1702 [math.RT]}}.
%Type = Article
\bibitem[{Kulish et~al.(2012)Kulish, Lyakhovsky, and
  Postnova}]{LyakhovskyPostnova2011}
\bibinfo{author}{P.~Kulish}, \bibinfo{author}{V.~Lyakhovsky},
  \bibinfo{author}{O.~Postnova},
\newblock \bibinfo{title}{Tensor powers for non-simply laced lie algebras
  b2-case},
\newblock \bibinfo{journal}{Journal of Physics: Conference Series}
  \bibinfo{volume}{346} (\bibinfo{year}{2012}) \bibinfo{pages}{012012}.
%Type = Article
\bibitem[{Kac and Wakimoto(1988)}]{kac1988modular}
\bibinfo{author}{V.~Kac}, \bibinfo{author}{M.~Wakimoto},
\newblock \bibinfo{title}{{Modular and conformal invariance constraints in
  representation theory of affine algebras}},
\newblock \bibinfo{journal}{Advances in mathematics(New York, NY. 1965)}
  \bibinfo{volume}{70} (\bibinfo{year}{1988}) \bibinfo{pages}{156--236}.
%Type = Article
\bibitem[{Walton(1989)}]{walton1989conformal}
\bibinfo{author}{M.~Walton},
\newblock \bibinfo{title}{{Conformal branching rules and modular invariants}},
\newblock \bibinfo{journal}{Nuclear Physics B} \bibinfo{volume}{322}
  (\bibinfo{year}{1989}) \bibinfo{pages}{775--790}.
%Type = Book
\bibitem[{Shifrin(2009)}]{shifrinmathematica}
\bibinfo{author}{L.~Shifrin}, \bibinfo{title}{Mathematica programming: an
  advanced introduction}, \bibinfo{year}{2009},
  \url{http://mathprogramming-intro.org/}.
%Type = Book
\bibitem[{Maeder(2000)}]{maeder2000computer}
\bibinfo{author}{R.~Maeder}, \bibinfo{title}{Computer science with
  Mathematica}, \bibinfo{publisher}{Cambridge University Press},
  \bibinfo{year}{2000}.
%Type = Inproceedings
\bibitem[{Casselman(1995)}]{casselman1995automata}
\bibinfo{author}{W.~Casselman},
\newblock \bibinfo{title}{Automata to perform basic calculations in coxeter
  groups},
\newblock in: \bibinfo{booktitle}{Representations of groups: Canadian
  Mathematical Society annual seminar, June 15-24, 1994, Banff, Alberta,
  Canada}, \bibinfo{organization}{Canadian Mathematical Society},
  \bibinfo{year}{1995}, volume~\bibinfo{volume}{16}, p.~\bibinfo{pages}{35}.
%Type = Article
\bibitem[{van~der Kallen(2011)}]{KallenShortlex}
\bibinfo{author}{W.~van~der Kallen},
\newblock \bibinfo{title}{Computing shortlex rewite rules with mathematica}
  (\bibinfo{year}{2011}),
  \url{http://www.staff.science.uu.nl/~kalle101/ickl/shortlex.html}.
%Type = Article
\bibitem[{Kazhdan and Lusztig(1994)}]{kazhdan1994tensor3}
\bibinfo{author}{D.~Kazhdan}, \bibinfo{author}{G.~Lusztig},
\newblock \bibinfo{title}{Tensor structures arising from affine lie algebras.
  iii},
\newblock \bibinfo{journal}{Journal of the American Mathematical Society}
  \bibinfo{volume}{7} (\bibinfo{year}{1994}).
%Type = Article
\bibitem[{Kazhdan and Lusztig(1993{\natexlab{a}})}]{kazhdan1993tensor1}
\bibinfo{author}{D.~Kazhdan}, \bibinfo{author}{G.~Lusztig},
\newblock \bibinfo{title}{Tensor structures arising from affine lie algebras.
  i},
\newblock \bibinfo{journal}{Journal of the American Mathematical Society}
  \bibinfo{volume}{6} (\bibinfo{year}{1993}{\natexlab{a}})
  \bibinfo{pages}{905--947}.
%Type = Article
\bibitem[{Kazhdan and Lusztig(1993{\natexlab{b}})}]{kazhdan1993tensor2}
\bibinfo{author}{D.~Kazhdan}, \bibinfo{author}{G.~Lusztig},
\newblock \bibinfo{title}{Tensor structures arising from affine lie algebras.
  ii},
\newblock \bibinfo{journal}{Journal of the American Mathematical Society}
  \bibinfo{volume}{6} (\bibinfo{year}{1993}{\natexlab{b}}).
%Type = Article
\bibitem[{Kac and Peterson(1984)}]{kac1984infinite}
\bibinfo{author}{V.~Kac}, \bibinfo{author}{D.~Peterson},
\newblock \bibinfo{title}{{Infinite-dimensional Lie algebras, theta functions
  and modular forms}},
\newblock \bibinfo{journal}{Adv. in Math} \bibinfo{volume}{53}
  (\bibinfo{year}{1984}) \bibinfo{pages}{125--264}.
%Type = Article
\bibitem[{Hwang and Rhedin(1995)}]{Hwang:1994yr}
\bibinfo{author}{S.~Hwang}, \bibinfo{author}{H.~Rhedin},
\newblock \bibinfo{title}{{General branching functions of affine Lie
  algebras}},
\newblock \bibinfo{journal}{Mod. Phys. Lett.} \bibinfo{volume}{A10}
  (\bibinfo{year}{1995}) \bibinfo{pages}{823--830},
  \href{http://arxiv.org/abs/hep-th/9408087}{{\tt arXiv:hep-th/9408087}}.
%Type = Article
\bibitem[{Hwang and Rhedin(1993)}]{hwang1993brst}
\bibinfo{author}{S.~Hwang}, \bibinfo{author}{H.~Rhedin},
\newblock \bibinfo{title}{The brst formulation of g/h wznw models},
\newblock \bibinfo{journal}{Nuclear Physics B} \bibinfo{volume}{406}
  (\bibinfo{year}{1993}) \bibinfo{pages}{165--184}.

\end{thebibliography}
\bibliographystyle{model1-num-names-arxiv-modified}






% \end{thebibliography}

\end{document}